\documentclass{article}
\usepackage{amsthm,amsmath,amssymb}
\usepackage{color}
\hoffset - 20 mm

\catcode`\@=11 \@addtoreset{equation}{section}

\catcode`\@=12

\newtheorem{Theorem}{Theorem}[section]
\newtheorem{Remark}{Remark}[section]
\newtheorem{Proposition}{Proposition}[section]
\newtheorem{Lemma}{Lemma}[section]
\newtheorem{Corollary}{Corollary}[section]

\newcommand{\bTheorem}[1]{
\begin{Theorem} \label{T#1} }
\newcommand{\eT}{\end{Theorem}}

\newcommand{\bRemark}[1]{
\begin{Remark} \label{T#1} }
\newcommand{\eR}{\end{Remark}}

\newcommand{\bProposition}[1]{
\begin{Proposition} \label{P#1}}
\newcommand{\eP}{\end{Proposition}}

\newcommand{\bLemma}[1]{
\begin{Lemma} \label{L#1} }
\newcommand{\eL}{\end{Lemma}}

\newcommand{\bCorollary}[1]{
\begin{Corollary} \label{C#1} }
\newcommand{\eC}{\end{Corollary}}

\newcommand{\bFormula}[1]{
\begin{equation} \label{#1}}
\newcommand{\eF}{\end{equation}}

\newcommand{\Ov}[1]{\overline{#1}}
\newcommand{\DC}{C^\infty_c}

\newcommand{\vr}{\varrho}

\newcommand{\vt}{\vartheta}
\newcommand{\vu}{\vc{u}}

\newcommand{\vc}[1]{{\bf #1}}
\newcommand{\Div}{{\rm div}_x}
\newcommand{\Grad}{\nabla_x}

\newcommand{\tn}[1]{\mbox {\F #1}}
\newcommand{\dx}{{\rm d} {x}}
\newcommand{\ds}{{\rm d} {s}}
\newcommand{\dr}{{\rm d} {r}}

\newcommand{\dt}{{\rm d} t }

\newcommand{\TT}{\Omega}

\newcommand{\intO}[1]{\int_{\TT} #1 \ \dx}

\newcommand{\TS}{R^{3 \times 3}_{{\rm sym},0} }

\newcommand{\ep}{\varepsilon}

\font\F=msbm10 scaled 1000
\newcommand{\R}{\mbox{\F R}}

\definecolor{grey}{rgb}{0.85,0.85,0.85}
\date{}

\long\def\greybox#1{%
    \newbox\contentbox%
    \newbox\bkgdbox%
    \setbox\contentbox\hbox to \hsize{%
        \vtop{
            \kern\columnsep
            \hbox to \hsize{%
                \kern\columnsep%
                \advance\hsize by -2\columnsep%
                \setlength{\textwidth}{\hsize}%
                \vbox{
                    \parskip=\baselineskip
                    \parindent=0bp
                    #1
                }%
                \kern\columnsep%
            }%
            \kern\columnsep%
        }%
    }%
    \setbox\bkgdbox\vbox{
        \color{grey}
        \hrule width  \wd\contentbox %
               height \ht\contentbox %
               depth  \dp\contentbox
        \color{black}
    }%
    \wd\bkgdbox=0bp%
    \vbox{\hbox to \hsize{\box\bkgdbox\box\contentbox}}%
    \vskip\baselineskip%
}

\newcommand{\HH}{\mathbb{H}}

\newcommand{\RR}{\mathbb{R}}

\newcommand{\QQ}{\mathbb{Q}}
\newcommand{\BB}{\mathbb{B}}

\newcommand{\SSS}{\mathbb{S}}

\newcommand{\YY}{\mathbb{Y}}

\newcommand{\pp}{\mathbf{p}}
\newcommand{\de}{\partial}

\newcommand{\calF}{{\mathcal F}}

\newcommand{\rhs}{right hand side}

\def\ocirc#1{\ifmmode\setbox0=\hbox{$#1$}\dimen0=\ht0 \advance\dimen0
  by1pt\rlap{\hbox to\wd0{\hss\raise\dimen0
  \hbox{\hskip.2em$\scriptscriptstyle\circ$}\hss}}#1\else {\accent"17 #1}\fi}

\date{}
\begin{document}

%%%%%%%%%%%%%%%%%%%%%%%%%%%%%%%%

%%%%%%%%%%%%%%%%%%%%%%%%%%%%%%%%

\title{Evolution of non-isothermal Landau-de Gennes\\ nematic liquid crystals flows with singular potential}
\author{Eduard Feireisl\thanks{Institute of Mathematics of the Czech
Academy of Sciences, \v Zitn\' a 25, 115 67 Praha 1, Czech
Republic. E-mail: {\tt feireisl@math.cas.cz} \ . The work of E.F.
was supported by Grant 201/09/0917 of GA \v CR. The work of E.F.
was partially supported by the FP7-IDEAS-ERC-StG \#256872
(EntroPhase).} \and Elisabetta Rocca\thanks{Dipartimento
di Matematica, Universit\`a di Milano, Via Sal\-di\-ni 50,
20133 Milano, Italy. E-mail {\tt elisabetta.rocca@unimi.it}\ . The
work of E.R. was supported by the FP7-IDEAS-ERC-StG \#256872
(EntroPhase).} \and Giulio Schimperna\thanks{Di\-par\-ti\-men\-to di
Ma\-te\-ma\-ti\-ca, U\-ni\-ver\-si\-t\`a 
di Pa\-vi\-a, Via Fer\-ra\-ta 1,
27100 Pavia,  Italy. E-mail: {\tt giusch04@unipv.it}\ . The work
of G.S. was supported by the MIUR-PRIN Grant 2008ZKHAHN ``Phase
transitions, hysteresis and multiscaling''.}  \and Arghir
Zarnescu\thanks{Pevensey III University of Sussex Falmer, BN1 9QH,
UK. E-mail: {\tt A.Zarnescu@sussex.ac.uk} \ . The work of A.Z. was
partially supported by the FP7-IDEAS-ERC-StG \#256872
(EntroPhase).}} \maketitle

\begin{abstract}
We discuss a $3D$ model describing the time evolution of nematic liquid crystals
in the framework of Landau-de Gennes theory, where the natural physical constraints
are enforced by a singular free energy bulk potential proposed by J.M. Ball and A. Majumdar.
The thermal effects are present through the component of the free energy that accounts for intermolecular
interactions. The model is consistent with the general principle of thermodynamics and mathematically tractable.
We identify the {\it a priori} estimates for the associated system of evolutionary partial
differential equations  and construct global-in-time weak solutions for arbitrary physically relevant initial data.
\end{abstract}

\medskip

\medskip

\tableofcontents

\section{Introduction}
\label{i}

The main aim of this paper is to derive and analyze
a thermodynamically consistent system of evolutionary equations
describing the dynamics of  nematic liquid crystal flows in $3D$.
We use the  abstract thermodynamic framework proposed by Fr\' emond \cite{Frem}
in conjunction with Beris-Edwards formulation \cite{BeEd} of  isothermal liquid crystal hydrodynamics (cf., e.g.,  \cite{DeToYe2}).
  The state of the complex fluid at the time $t$ and  spatial position
$x$ is described by means of a \emph{$Q$-tensor} field $\tn{Q} = \tn{Q}(t,x)$ for the nematic director orientation,
the \emph{velocity field} $\vu = \vu(t,x)$, and the \emph{absolute temperature}
$\vt = \vt(t,x)$.

The main characteristic of nematic liquid crystals is the locally preferred
orientation of the nematic molecule directors. This can be  described  by the
$Q$-tensors, that are suitably normalized second order moments of the probability
distribution function of the molecules. More precisely if $\mu_x$ is a probability
measure on the unit sphere $\mathbb{S}^2$, representing the orientation of the molecules at a point
$x$ in space, then a $Q$-tensor $\tn{Q}(x)$ is a symmetric and traceless $3\times 3$ matrix defined as
\begin{equation}\nonumber
  \tn{Q}(x)=\int_{\mathbb{S}^2}\left(\pp\otimes \pp-\frac{1}{3}\tn{I}\right)\,d\mu_x(\pp)
\end{equation}
and it is supposed to be a  a crude measure of how the probability measure $\mu_x$ deviates
from the isotropic measure $\bar\mu$ where $d\bar\mu=\frac{1}{4\pi}dA$, see \cite{DeGennes}.
In the Onsager model (cf. \cite{DeGennes}, \cite{MaSa})
the probability measure is assumed to be continuous with density $\rho=\rho(\pp)$. Then
\begin{equation}\label{Qrep}
  \tn{Q}(x)=\int_{\mathbb{S}^2}\left(\pp\otimes \pp-\frac{1}{3}\tn{I}\right)\,\rho(\pp)\, d\pp\,.
\end{equation}
The fact that $\mu_x$ is a probability measure imposes a constraint on the eigenvalues of
$\tn{Q}$, namely that they are bound between the values $-\frac{1}{3}$ and $\frac{2}{3}$,
see \cite{BalMaj}. Thus not any traceless $3\times 3$ matrix is a {\it physical} $Q$-tensor
but only those whose eigenvalues  are in $(-\frac{1}{3},\frac{2}{3})$. The hydrodynamic models available
in the literature do not have, to our knowledge, a natural way of preserving this physical
eigenvalue constraint on the traceless and symmetric matrices. One possibility is to use a
singular potential, such as the one proposed by  Ball and Majumdar \cite{BalMaj},
that enforces the  physical constraints satisfied by the $Q$-tensors, and this is the
solution we adopt in our model.

The hydrodynamic theory for the $Q$-tensorial isothermal model in case of regular bulk  potential
$$\psi_B(\tn{Q})=\frac{a}{2}{\rm tr}(\tn{Q}^2)-\frac{b}{2}{\rm tr}(\tn{Q}^3)+\frac{c}{4}{\rm tr}^2(\tn{Q}^2)$$
has been recently studied in \cite{PaZa2011} and \cite{PaZa2012}.

Let us notice that in the literature there are few papers dealing with non-isothermal models for
liquid crystal dynamics. Two attempts, in the case when the evolution of the director is
described by the vectorial director field $\vc{d}$ (standing for the preferred orientation
of the molecules at any point), were made in~\cite{FRS} and~\cite{FFRS}. In particular,
in~\cite{FFRS} the stretching and rotation effects of the director field induced by
the straining of the fluid were considered and the existence of global in time weak
solutions was obtained for the corresponding initial boundary value problem.
In the present contribution we  follow the thermodynamic approach exploited in~\cite{FFRS}
in order to deal with the tensorial model obtained using a non-isothermal version of the
singular bulk potential $f(\tn{Q})$ proposed in~\cite{BalMaj}. In the spirit of \cite{DeGennes},
we include the temperature dependence in the potential assuming that the coupling term
in the free energy functional is given by (cf. also \cite{SchoSlu} and  \cite{MottramNewton})
\bFormula{psi}
  \psi_B (\vt, \tn{Q}) = f(\tn{Q}) - U(\vt) G(\tn{Q}).
\eF
Hence, $\psi_B$ is the sum of a {\it singular}\ part $f$, independent of temperature $\vt$,
with a smooth perturbation depending both on $\vt$ and on $\tn{Q}$. We assume $U$ to be a
convex and decreasing function of $\vt$ having controlled growth at infinity
(cf.~(\ref{i3}--\ref{i3bis}) below). Actually, in a neighbourhood of a characteristic temperature
$\vt^*$ of the crystal the function $U$ can display a linear growth and can change sign at $\vt^*$, like
$U(\vt)=\alpha(\vt^*-\vt)$. According to \cite{DeGennes}, the function $G$ can be, e.g., given by
$G(\tn{Q})= \tn{Q}_{ij}\tn{Q}_{ij}=\text{tr}(\tn{Q}^2)$. Let us notice that in the Ball and Majumdar
\cite{BalMaj} approach the temperature dependence is of a different type: actually, they assume
\bFormula{psimaj}
  \psi_B (\vt, \tn{Q}) = \vt f(\tn{Q})+ G(\tn{Q}).
\eF
However in \cite{BalMaj} only the stationary case is considered; hence, one can freely
divide by $\vt$ the expression $\psi_B$ in \eqref{psimaj} and obtain an expression which corresponds
to \eqref{psi} at least for values of $\vt$ not too distant from the critical
temperature $\vt^*$. It is worth noting that, in the evolutionary setting, dealing with
a free energy of the form \eqref{psimaj} would be mathematically more complicated since the
coupling occurs in the singular part of the potential; for this reason we expect that
weaker analytic results could be proved in that case. We will devote a forthcoming
paper to the analysis of the evolutive model with the free energy \eqref{psimaj}.

Comparing the present analysis with the previous paper \cite{FFRS},
a major difficulty is provided here by the presence of the singular potential $f$,
which has to be properly handled by means of convex analysis tools.
In addition to that, we consider here a more complicated version of
the heat equation involving an explicit dependence of the thermal
relaxation coefficient with respect to the $Q$-tensor (cf.~\eqref{p6} below).
Actually, this choice, which is more realistic from the physical point of view
(in particular, it gives rise to an entropy $s$ depending also
on $\tn{Q}$ and not only on $\vt$, cf.~\eqref{i27} below), creates
a number of additional mathematical difficulties. The key
point, which requires some care to be accomplished, is related to the proof
of \emph{regularity}\/ and \emph{strict positivity}\/ of $\vt$ at the
approximate level, two properties which are crucial for the purposes
of proving the validity of the entropy inequality and
of the total energy balance in the frame of weak solutions.

%%%%%%%%%%%%%%%%%%%%%%%%%%%%%%%%%%%%%%%%%%%%%%%%%%%%%%%%%%%%%%%%%%%%%%%%%%%%%%%%%%%%%%%%%%%%%%%%%%%%%%%%%

\subsection{Landau-de Gennes free energy with the Ball-Majumdar bulk potential}

Denote by $\TS$ the linear space of symmetric traceless $3\times 3$ real-valued matrices.
The \emph{Landau-de Gennes free energy} takes the form
\bFormula{i1}
  \mathcal{F}(\tn{Q}, \Grad \tn{Q}, \vt ) = \frac{1}{2} |\Grad \tn{Q} |^2 + \psi_B (\vt, \tn{Q}) - \vt \log(\vt) ,
\eF
where $\tn{Q}(x)\in \TS$ for all $x$ in the smooth domain $\Omega\subset\R^3$.

%\[
%\mathcal{L} [\tn{Q}] = \tn{Q} - \frac{1}{3} {\rm Tr}[\tn{Q}] \tn{I}, \ \tn{Q} \in R^{3 \times 3}_{\rm sym},
%\]
%denotes the projection onto the space of traceless tensors.

Ball and Majumdar \cite{BalMaj}
introduced the bulk component of the internal energy functional by means of a
singular functional $\psi_B = \psi_B(\vt, \tn{Q})$ that, for any fixed temperature $\vt$, blows up
when at least one of the eigenvalues of $\tn{Q}$ approaches the limiting value
$- 1/3$. In particular, the boundedness of the free energy enforces the boundedness of $\tn{Q}$ in $L^\infty$. Specifically, we set
\bFormula{i2}
\psi_B (\vt, \tn{Q}) = f(\tn{Q}) - U(\vt) G(\tn{Q}) \ \mbox{for}\ \tn{Q} \in \TS,
\eF
where
\[
f(\tn{Q}) = \left\{ \begin{array}{l} \inf_{\rho \in \mathcal{A}_{\tn{Q}}}
\int_{S^2} \rho(\vc{p}) \log (\rho(\vc{p})) \ {\rm d}\vc{p} \ \mbox{if} \
\lambda_i[\tn{Q}] \in (-1/3,2/3),\ i=1,2,3,
\\ \\ \infty \ \mbox{otherwise,}\end{array} \right.
\]
\[
\mathcal{A}_{\tn{Q}} = \left\{ \rho: S^2 \to [0, \infty) \ \Big|\
\rho\in L^1(S^2),\ \int_{S^2} \rho (\vc{p}) \ {\rm d} \vc{p} = 1;
\tn{Q} = \int_{S^2} \left( \vc{p} \otimes \vc{p} - \frac{1}{3}
\tn{I} \right) \rho(\vc{p}) \ {\rm d} \vc{p} \right\}.
\]

The function $f$ is the singular component of the bulk potential. In here singular refers to the fact that
the domain is not the whole space (while inside the domain the function is in fact smooth). The function $f$
enjoys the following properties that can easily be deduced  from \cite[Section 3, Prop. 1]{BalMaj}:
\begin{itemize}
\item $f : \TS \to [-K, \infty]$ is convex and lower semi-continuous, with $K \ge 0$.
\item The domain of $f$,
\[
\mathcal{D}[f] = \{ \tn{Q} \in R^{3 \times 3}_{\rm sym,0} \ | \ f(\tn{Q}) < \infty \}
= \{  \tn{Q} \in R^{3 \times 3}_{\rm sym,0} \ | \ \lambda_i[\tn{Q}] \in (- 1/3,2/3) \},
\]
is an open convex subset of $\TS$.
\item $f$ is smooth in $\mathcal{D}[f]$.
\end{itemize}

The potential $G$ characterizes the action of intermolecular forces. In contrast with Ball and Majumdar \cite{BalMaj},
we suppose the temperature changes act on this component of the bulk potential. Such a hypothesis is quite common
in the literature (see for instance \cite{DeGennes}); actually, $U$ is typically assumed to change sign at a
critical temperature. Here, we assume that there exists a positive constant $c$ such that
\bFormula{i3}
U \in C^0[0, +\infty)\cap C^2(0,+\infty), \ U(0) > 0 , \ U' \leq 0, \ U \ \mbox{convex in}\ [0, \infty),
 \ \limsup_{ \vt \to \infty } U''(\vt) \vt^{3/2} < + \infty,
\eF
\bFormula{i3bis}
 | U'(\vt) | \le c |\vt|^{-1/2} \ \hbox{ for all } \vt \in(0,\infty),
\eF
and
\bFormula{i4}
  G \in C^3 (\TS), \ G \geq 0, \  G(\tn{Q})= G(R\tn{Q}R^t) \  \hbox{for all } R\in SO(3).
\eF
Let us note that the choice $G(\tn{Q})= \text{tr}(\tn{Q}^2)$ corresponds to that of \cite{DeGennes}.
Moreover, it is not restricive to assume that $G$ is uniformly bounded together with his first and
second derivatives.

\subsection{Thermodynamics}

In accordance with the general thermodynamical framework of  \cite{Frem}, we introduce
the set of \emph{state variables}
\[
E = \left( \tn{Q} , \Grad \tn{Q} , \vt \right),
\]
together with the \emph{dissipative variables}
\[
  E^{\rm d} = \left( \ep ( \vu ) , \frac{D \tn{Q}}{D t}, \Grad \vt \right),
\]
where
\[
\ep (\vu) = \frac{1}{2} \left( \Grad \vu + \Grad^t \vu \right)
\]
is the symmetric velocity gradient, and
\begin{equation}\label{dS}
\frac{D \tn{Q} }{Dt} \equiv \partial_t \tn{Q} + \vu \cdot \Grad \tn{Q} -
\tn{S} ( \Grad \vu, \tn{Q} )
\end{equation}
is an analogue of \emph{material derivative}
characterizing the time evolution of the tensor $\tn{Q}$,
with
\bFormula{i5}
\tn{S}(\Grad \vu, \tn{Q}) = \left( \xi \ep(\vu) + \omega (\vu) \right)\left(
\tn{Q} + \frac{1}{3} \tn{I} \right) + \left( \tn{Q} + \frac{1}{3} \tn{I} \right) \left( \xi \ep (\vu) - \omega (\vu) \right) -
2 \xi \left( \tn{Q} + \frac{1}{3} \tn{I} \right) \left( \QQ: \Grad \vu \right),
\eF
\[
\omega(\vu) = \frac{1}{2} \left( \Grad \vu - \Grad^t \vu \right),
\]
where $\xi$ is a fixed scalar parameter, measuring the ratio between the rotation and the aligning effect that a shear flow would have over the directors, see Beris and Edwards \cite{BeEd}.

The evolution is ruled by the
\emph{pseudopotential of dissipation} $\Phi$
\bFormula{i6}
\Phi (E^{\rm d},E) = \mu(\vt) |\ep (\vu)|^2 + I_0 (\Div
\vu) + \frac{\kappa(\vt)}{2\vt} |\Grad \vt |^2 + \frac{1}{2 \Gamma
(\vt)} \left| \frac{ D \tn{Q} }{Dt} \right|^2,
\eF
with the shear viscosity coefficient $\mu$, the heat conductivity
coefficient $\kappa$, and the collective rotational viscosity
coefficient $\Gamma$. The \emph{incompressibility} of the fluid is
formally enforced by the functional $I_0$ - the indicator function
of the point $\{ 0 \}$,
\[
I_0 (z) = \begin{cases} 0 & \ \mbox{if}\ z = 0, \\
+ \infty & \ \mbox{otherwise}. \end{cases}
\]

Similarly to $U$, the transport coefficients $\mu(\vt)$, $\kappa(\vt)$, and $\Gamma(\vt)$ change with
temperature. For the sake of simplicity, we suppose that
\bFormula{i7}
\mu, \ \kappa, \ \Gamma \in C^2[0,\infty), \ \left\{
\begin{array}{c} 0 < \underline{\mu} \leq \mu(\vt) \leq \Ov{\mu}, \\ \\
0 < \underline{\kappa} \leq \kappa(\vt) \leq \Ov{\kappa}, \\ \\
0 < \underline{\Gamma} \leq \Gamma(\vt) \leq \Ov{\Gamma}
\end{array} \right\} \ \mbox{for all}\ \vt \geq 0.
\eF

%%%%%%%%%%%%%%%%%%%%%%%%%%%%%%%%%%%%%%%%%%%%%%%%%%%%%%%%%%%%%%%%%%%%%%%%%%%%%%%%%%%%%%%%%%%%%%%%%%%%%%

\subsection{Time evolution}

We assume that the fluid has a constant density ${\vr}$, say ${\vr} = 1$.
Then, in accordance with the general principles developed
in Fr{\' e}mond \cite[Chapters 2,3]{Frem}, the time evolution of the
system is uniquely determined by the choice of the potentials $\mathcal{F}$ and $\Phi$.

\subsubsection{Momentum equation}

Newton's second law is expressed by means of a modified
\emph{Navier-Stokes system}:
\bFormula{i8}
\partial_t \vu + \Div (\vu \otimes \vu) =
\Div \sigma + \vc{g},
\eF
where $\vc{g}$ is a driving force, and $\sigma$ denotes the stress tensor,
\[
\sigma = \sigma^{\rm d} + \sigma^{\rm nd}.
\]
The dissipative component of the stress reads
\[
\sigma^{\rm d} = \frac{\partial \Phi}{\partial \ep(\vu)} =
\frac{ \mu(\vt) }{2}
(\Grad \vu + \Grad^t \vu) - p \tn{I},
\]
where $p$ is the pressure. We have
\[
- p \in \partial I_0 (\Div \vu)
\]
yielding the standard incompressibility constraint
\bFormula{i9}
\Div \vu = 0.
\eF
The specific form of $\sigma^{\rm nd}$ will be derived below.

%%%%%%%%%%%%%%%%%%%%%%%%%%%%%%%%%%%%%%%%%%%%%%%%%%%%%%%%%%%%%%%%%%%%%%%%%%%%%%%%%%%%%%%%%%

\subsubsection{Entropy production}

The heat flux $\vc{q}$ can be decomposed as
\[
\vc{q} = \vc{q}^{\rm d} + \vc{q}^{\rm nd},
\]
where the dissipative component obeys the standard \emph{Fourier law}
\[
\vc{q}^{\rm d} = - \vt \frac{\partial \Phi}{\partial \Grad \vt} = -
\kappa(\vt) \Grad \vt,
\]
with the associated entropy flux $\vc{q}_e = \vc{q}^{\rm d}/ \vt$.

The energy density tensor is taken to be
\[
\tn{B} = \tn{B}^{\rm d} + \tn{B}^{\rm nd},
\]
\bFormula{i10}
\tn{B}^{\rm d} = \frac{\partial \Phi}{\partial \frac{D \tn{Q}}{Dt} }
= \frac{1}{\Gamma (\vt)} \frac{D \tn{Q}}{Dt}, \quad
\tn{B}^{\rm nd} =
\frac{ \partial \mathcal{F} }{\partial \tn{Q}} = \mathcal{L} \left[
\frac{ \partial f (\tn{Q}) }{\partial \tn{Q} } \right] - U(\vt) \mathcal{L} \left[
\frac{ \partial G (\tn{Q}) }{\partial \tn{Q} } \right],
\eF
where
\[
\mathcal{L} [h(\tn{Q})] = h(\tn{Q}) - \frac{1}{3} {\rm tr}[h(\tn{Q})] \tn{I},\text{ for any }h(\tn{Q}) \in R^{3 \times 3}_{\rm sym},
\]
denotes the projection onto the space of traceless tensors. In other words, $\tn{B}^{\rm nd}$ can
be seen as the subdifferential of $\calF$ with respect to $\QQ$ in that space.

We state the \emph{entropy equation} in the form (cf. \cite{bcf, bf} for a complete derivation of this equation)
\bFormula{i11}
 \partial_t s + \Div (s \vu)  - \Div \vc{q}_e = \frac{1}{\vt} \left(
\sigma^{\rm d} : \ep(\vu) + \tn{B}^{\rm d} : \frac{D \tn{Q}}{Dt} +
\frac{\kappa(\vt)}{\vt} |\Grad \vt |^2 \right)
\eF
\[
= \frac{1}{\vt} \left( \mu(\vt) |\ep(\vu)|^2 +
\frac{1}{\Gamma(\vt)} \left| \frac{ D \tn{Q} }{Dt} \right|^2
+\frac{\kappa(\vt)}{\vt} |\Grad \vt |^2 \right) \geq 0,
\]
with the \emph{entropy}
\bFormula{entro}
  s = - \frac{\partial \mathcal{F}}{\partial \vt} = 1 + \log(\vt) +  U'(\vt) G(\tn{Q}).
\eF
Note that, in accordance with hypotheses (\ref{i3}), (\ref{i4}), the entropy $s$
is an increasing function of the
temperature.

%%%%%%%%%%%%%%%%%%%%%%%%%%%%%%%%%%%%%%%%%%%%%%%%%%%%%%%%%%%%%%%%%%%%%%%%%%%%%%%

\subsubsection{$Q$-tensor evolution}

The internal energy balance reads
\bFormula{i12}
 \partial_t e + \Div (e \vu)  + \Div \vc{q} =
\sigma : \Grad \vu + \tn{B}: \frac{D \tn{Q}}{Dt} + \tn{Y} : \Grad
\frac{D \tn{Q}}{Dt},
\eF
where
\bFormula{i13}
 e =  \mathcal{F} +  \vt s = \frac{1}{2} |\Grad \tn{Q}|^2 + f(\tn{Q})-
\Big( U(\vt) - \vt U'(\vt) \Big)
G(\tn{Q}) + \vt,
\eF

\noindent
and $\tn{Y}$ is the energy flux tensor
\bFormula{i14}
\tn{Y} = \tn{Y}^{\rm nd} = \frac{\partial \mathcal{F} }{\partial
\Grad \tn{Q}} = \Grad \tn{Q}.
\eF
For simplicity, we assume here that $\tn{Y}$ has no dissipative component
(and, correspondingly, that $\Phi$ is independent of $\Grad \tn{Q}_t$).

The principle of virtual powers (see Fr{\' e}mond
\cite[Chapter 2]{Frem}) yields the time evolution of $\tn{Q}$, namely,
\bFormula{i15}
\Div \tn{Y} = \tn{B};
\eF
in other words,
\bFormula{i16}
\partial_t \tn{Q} + \Div (\tn{Q} \vu) - \tn{S}(\Grad \vu, \tn{Q}) =
\Gamma(\vt) \left( \Delta \tn{Q} - \mathcal{L} \left[ \frac{\partial
f (\tn{Q})}{\partial \tn{Q}} \right] + U(\vt)\mathcal{L} \left[ \frac{\partial
G (\tn{Q})}{\partial \tn{Q}} \right]  \right),
\eF
where $\tn{S}$ is defined in \eqref{i5}.
It is easy to check that the space $\TS$ is invariant for solutions of
(\ref{i16}), specifically $\tn{Q}(t, \cdot) \in \TS$ for any $t \geq 0$ as soon as
$\tn{Q}(0, \cdot) \in \TS$.

%%%%%%%%%%%%%%%%%%%%%%%%%%%%%%%%%%%%%%%%%%%%%%%%%%%%%%%%%%%%%%%%%%%%%%%%%%%%%%%%%%%%%%

\subsubsection{Total energy balance}

Taking the scalar product of the momentum equation (\ref{i8}) with
$\vu$ and adding the resulting expression to (\ref{i12}), we obtain
the total energy balance in the form
\bFormula{i17}
\partial_t \left( \frac{1}{2} |\vu |^2 + e \right) + \Div
\left(\left(\frac{1}{2} |\vu |^2 + e\right) \vu \right) + \Div \vc{q}
\eF
\[
= \Div (\sigma \vu) + \Div \left( \Gamma (\vt) \Grad \tn{Q} :
\left( \Delta \tn{Q} - \mathcal{L} \left[ \frac{\partial f
(\tn{Q})}{\partial \tn{Q}} \right] + U(\vt) \mathcal{L} \left[ \frac{\partial G
(\tn{Q})}{\partial \tn{Q}} \right]  \right)   \right) +
\vc{g} \cdot \vu.
\]
It remains to determine $\vc{q}^{\rm nd}$ and $\sigma^{\rm nd}$.

To this end, we first multiply the entropy balance \eqref{i11}
by $\vt$. This gives
\bFormula{g11}
  \vt s_t - s \vu \cdot \Grad \vt+ \Div (\vt s \vu)
       +\Div {\bf q}^d
    = \sigma^{\rm d} : \ep(\vu) + \BB^{\rm d} : \frac{D \tn{Q}}{Dt}.
\eF
Next, using \eqref{entro}, \eqref{i13}, \eqref{i10} and \eqref{i14},
we get
\bFormula{g12}
  e_t = \YY :: \de_t \Grad \QQ + \BB^{\rm nd} : \de_t \QQ + \vt s_t.
\eF
Moreover,
\bFormula{g13}
  \Div ( e \vu ) = \Div \big( ( \calF + s \vt ) \vu \big)
   = \vu \cdot \big( \calF_{\vt} \Grad \vt + \calF_{\QQ} : \Grad \QQ
     + \calF_{\Grad \QQ} :: \Grad \Grad \QQ \big)
   + \Div ( s \vt \vu )
\eF
\[
  = - s \vu \cdot \Grad \vt
   + \vu \cdot ( \BB^{\rm nd} : \Grad \QQ )
    + \vu \cdot ( \YY :: \Grad \Grad \QQ \big)
  + \Div ( s \vt \vu ).
\]
Replacing (\ref{g12}--\ref{g13}) in \eqref{i12} and using \eqref{dS},
we then have
\bFormula{g14}
   \YY :: \de_t \Grad \QQ + \BB^{\rm nd} : \de_t \QQ + \vt s_t
   - s \vu \cdot \Grad \vt
   + \vu \cdot ( \BB^{\rm nd} : \Grad \QQ )
   + \vu \cdot ( \YY :: \Grad \Grad \QQ \big)
  + \Div ( s \vt \vu )
  + \Div \vc{q}
\eF
\[
 = \sigma^{\rm d} : \Grad \vu
 + \sigma^{\rm nd} : \Grad \vu
 + \tn{B}^{\rm d}: \frac{D \tn{Q}}{Dt}
 + \tn{B}^{\rm nd}: \frac{D \tn{Q}}{Dt}
 + \tn{Y} :: \Grad \frac{D \tn{Q}}{Dt}
\]
\[
 = \sigma^{\rm d} : \Grad \vu
 + \sigma^{\rm nd} : \Grad \vu
 + \tn{B}^{\rm d}: \frac{D \tn{Q}}{Dt}
\]
\[
 \mbox{}
 + \tn{B}^{\rm nd}: \big(\de_t \QQ + \vu \cdot \Grad \QQ - \SSS(\Grad \vu, \QQ) \big)
 + \tn{Y} :: \big( \de_t \Grad \QQ + \Grad ( \vu \cdot \Grad \QQ )
  - \Grad  \SSS(\Grad \vu, \QQ) \big).
\]
Simplifying some terms and using symmetry of $\sigma^{\rm d}$, we have
more precisely
\bFormula{g14b}
   \vt s_t
   -   \vu s  \cdot\Grad \vt
  + \Div ( s \vt \vu )
  + \Div \vc{q}^{\rm d}
  + \Div \vc{q}^{\rm nd}
\eF
\[
 = \sigma^{\rm d} : \ep(\vu)
 + \sigma^{\rm nd} : \Grad \vu
 + \tn{B}^{\rm d}: \frac{D \tn{Q}}{Dt}
 - \tn{B}^{\rm nd}: \SSS(\Grad \vu, \QQ)
 + \tn{Y} :: (\Grad \vu \cdot \Grad \QQ )
 - \tn{Y} :: \Grad  \SSS(\Grad \vu, \QQ).
\]
Then, subtracting \eqref{g11} from \eqref{g14b} and using \eqref{i10},
we arrive at
\bFormula{es8}
  \Div \vc{q}^{\rm nd}
\eF
\[
  = \sigma^{\rm nd} : \Grad \vu
   - \left(  \mathcal{L} \left[ \frac{\partial f(\tn{Q})}{\partial \tn{Q} } \right]
        - U(\vt) \mathcal{L} \left[ \frac{\partial G(\tn{Q})}{\partial \tn{Q} } \right] \right)
    : \tn{S} (\Grad \vu , \tn{Q} )
\]
\[
  - \Grad \tn{Q} :: \Grad \tn{S} (\Grad \vu, \tn{Q})
   + \left( \Grad \tn{Q} \odot \Grad \tn{Q} \right) : \Grad \vu
\]
\[
  = \sigma^{\rm nd} : \Grad \vu- \Div \left( \Grad \tn{Q} : \tn{S}(\Grad \vu , \tn{Q}) \right)
\]
\[
  + \left( \Delta \tn{Q} - \mathcal{L} \left[ \frac{\partial f(\tn{Q})}{\partial \tn{Q} } \right]
   + U(\vt) \mathcal{L} \left[ \frac{\partial G(\tn{Q})}{\partial \tn{Q} } \right] \right):
\tn{S} (\Grad \vu , \tn{Q} ) + \left( \Grad \tn{Q} \odot \Grad
\tn{Q} \right) : \Grad \vu.
\]
Consequently, we deduce that
\bFormula{i18}
  \vc{q}^{\rm nd} = - \Grad \tn{Q} : \tn{S}(\Grad \vu , \tn{Q}),
\eF
and
\bFormula{i19}
  \sigma^{\rm nd} = \tn{Q} \tn{H} - \tn{H} \tn{Q}
   + 2 \xi \left[ \tn{H}  : \tn{Q} \right] \left( \tn{Q} + \frac{1}{3} \tn{I} \right)
   - \xi \left[ \tn{H}  \left( \tn{Q} + \frac{1}{3} \tn{I} \right)
   + \left( \tn{Q} + \frac{1}{3} \tn{I} \right) \tn{H} \right]- \left( \Grad \tn{Q} \odot \Grad \tn{Q} \right),
\eF
where we have denoted
\bFormula{i20}
  \tn{H} \equiv \Delta \tn{Q}
  - \mathcal{L} \left[ \frac{\partial f(\tn{Q})}{\partial \tn{Q} } \right]
   + U(\vt) \mathcal{L} \left[ \frac{\partial G(\tn{Q})}{\partial \tn{Q} } \right]
\eF
and we have used the identity
\bFormula{mat}
 - \tn{H} : \tn{S}(\Grad \vu, \tn{Q})
\eF
\[
 = \left( \tn{Q} \tn{H} - \tn{H} \tn{Q} \right) : \Grad \vu + 2 \xi \left( \tn{H}: \tn{Q} \right) \left( \tn{Q} : \Grad \vu \right)
  - \xi \left[ \tn{H} \left( \tn{Q}  + \frac{1}{3} \tn{I} \right)
   + \left( \tn{Q}  + \frac{1}{3} \tn{I} \right) \tn{H} \right]: \Grad \vu
\]
that holds for any symmetric matrix $\tn{H}$.
%

%%%%%%%%%%%%%%%%%%%%%%%%%%%%%%%%%%%%%%%%%%%%%%%%%%%%%%%%%%%%%%%%%%%%%%%%%%%%%%%%%%%%%%

\subsubsection{Evolutionary system}

The computations given in the previous section
permit to write the resulting evolutionary system in a concise form:

\greybox{
\centerline{\textsc{Incompressibility:}}

\vspace{-7mm}

\bFormula{i21}
\Div \vu = 0;
\eF

\medskip

\centerline{\textsc{Momentum equation:}}

\vspace{-7mm}

\bFormula{i22}
\partial_t \vu + \Div (\vu \otimes \vu) =
\Div \sigma + \vc{g};
\eF

\medskip

\centerline{\textsc{Order parameter evolution:}}

\vspace{-7mm}

\bFormula{i23}
\partial_t \tn{Q} + \Div (\tn{Q} \vu ) - \tn{S} (\Grad \vu, \tn{Q}) =
\Gamma(\vt) \tn{H};
\eF

\medskip

\centerline{\textsc{Total energy balance:}}

\vspace{-7mm}

\bFormula{i24}
 \partial_t \left( \frac{1}{2} |\vu|^2 + e \right) +
\Div \left( \left( \frac{1}{2} |\vu|^2 + e \right) \vu \right)
+ \Div \vc{q}
\eF
\[
=\Div (\sigma \vu ) + \Div \Big( \Gamma (\vt) \Grad \tn{Q} : \tn{H}
  \Big) +
\vc{g} \cdot \vu,
\]
}

\medskip

\noindent
with the stress tensor
\bFormula{i25}
\sigma = \mu(\vt) \left( \Grad \vu + \Grad^t \vu \right) - p \tn{I}
\eF
\[
+ 2 \xi \left( \tn{H}: \tn{Q} \right) \left( \tn{Q} + \frac{1}{3} \tn{I} \right)
- \xi \left[ \tn{H} \left( \tn{Q}  + \frac{1}{3} \tn{I} \right) +
\left( \tn{Q}  + \frac{1}{3} \tn{I} \right) \tn{H} \right] +
\left( \tn{Q} \tn{H} - \tn{H} \tn{Q} \right) - \Grad \tn{Q} \odot \Grad \tn{Q},
\]
where
\[
  \tn{H} \equiv \Delta \tn{Q}
  - \mathcal{L} \left[ \frac{\partial f(\tn{Q})}{\partial \tn{Q} } \right]
   + U(\vt) \mathcal{L} \left[ \frac{\partial G(\tn{Q})}{\partial \tn{Q} } \right],
\]
and the internal energy flux
\bFormula{i26}
\vc{q} = - \kappa (\vt) \Grad \vt - \Grad \tn{Q} : \tn{S} (\Grad \vu, \tn{Q}),
\eF
where
\[
\tn{S}(\Grad \vu, \tn{Q}) = \left( \xi \ep(\vu) + \omega (\vu) \right)\left(
\tn{Q} + \frac{1}{3} \tn{I} \right) + \left( \tn{Q} + \frac{1}{3} \tn{I} \right) \left( \xi \ep (\vu) - \omega (\vu) \right) -
2 \xi \left( \tn{Q} + \frac{1}{3} \tn{I} \right) \left( \QQ: \Grad \vu \right)\,.
\]

We recall that
\bFormula{i27}
e = \frac{1}{2} |\Grad \tn{Q}|^2 + f(\tn{Q}) - \Big( U(\vt) - \vt U'(\vt) \Big)
G(\tn{Q}) + \vt, \ s =1 + \log(\vt) +  U'(\vt) G(\tn{Q}),
\eF
where we have anticipated that the relation
$\mathcal{L}[\tn{Q}] = \tn{Q}$ is preserved in the course of evolution.

The system (\ref{i21}--\ref{i24}) may be supplemented by the \emph{entropy inequality}
\bFormula{i28}
 \partial_t s + \Div (s \vu)  - \Div \left(
\frac{\kappa(\vt)}{\vt} \Grad \vt \right)
\eF
\[
\geq \frac{1}{\vt} \left( \frac{\mu(\vt)}{2} \left| \Grad \vu + \Grad^t \vu \right|^2 +
{\Gamma(\vt)} \left| \tn{H} \right|^2
+\frac{\kappa(\vt)}{\vt} |\Grad \vt |^2 \right).
\]

In order to avoid problems related to the presence of a kinematic boundary, we suppose that the fluid motion is
spatially periodic. This can be conveniently formulated by taking the spatial domain $\Omega\subset\R^3$ as a flat torus

\greybox{
\bFormula{i29}
\Omega = \left( [-\pi, \pi]|_{\{ -\pi, \pi \} } \right)^3.
\eF
}
\noindent
Note also that the pressure $p$ appears \emph{explicitly} in the energy balance (\ref{i24}), in particular, $p$ must be determined
from the Navier-Stokes system (\ref{i22}) by means of the Helmholtz projection. Such a step may involve insurmountable difficulties
in the case of general boundary conditions.

The original state of the system is given by the initial conditions

\greybox{
\bFormula{i30}
\vu(0, \cdot) = \vu_0, \ \tn{Q}(0, \cdot) = \tn{Q}_0, \ \vartheta(0, \cdot) = \vt_0.
\eF
}

Our goal in this paper is to study the initial-value problem (\ref{i21}--\ref{i24}), supplemented with
the boundary conditions (\ref{i29}) and the initial conditions (\ref{i30}), in the framework of weak solutions.
In the sequel we will show that, for any choice of finite energy initial data (cf.~Section~\ref{m} below),
the problem possesses a global-in-time weak solution, which additionally satisfies
the entropy inequality \eqref{i28} in the sense of distributions.
To this end, we first derive formal {\it a priori} bounds in order to
facilitate the reading of the main rather technical part of the proof, see Section \ref{a}.
The global-in-time weak solutions are constructed
as a limit of solutions of a family of approximate problems introduced in Section \ref{p}.
The most delicate part of the proof is showing strict positivity of the (absolute) temperature by means
of a weak variant of the parabolic comparison theorem. The proof of convergence of approximate solutions is completed in
Section~\ref{c}.

%%%%%%%%%%%%%%%%%%%%%%%%%%%%%%%%%%%%%%%%%%%%%%%%%%%%%%%%%%%%%%%%%%%%%%%%%%%%%%%%%%%%%%%%%%%%%%%%%%%%%%%%%%

\section{Weak solutions, main results}
\label{m}

Weak solutions to the problem (\ref{i21}--\ref{i24}), (\ref{i29}), (\ref{i30}) belong to the regularity classes indicated by the {\it a priori}
bounds discussed in Section \ref{a} below. In particular, we have
\bFormula{m1}
\vu \in L^\infty(0,T; L^2(\Omega; \R^3)), \ \left\{
\begin{array}{c} \tn{Q} \in L^\infty((0,T) \times \Omega; \TS), \\ \\ \Grad \tn{Q} \in L^\infty(0,T;
L^2(\Omega; R^{27})), \\ \\ f(\tn{Q}) \in L^\infty(0,T; L^1(\Omega)) \end{array} \right\}, \
\left\{ \begin{array}{c} \vt \in L^\infty(0,T; L^1(\Omega)), \\ \\
 \vt \in C_w([0,T]; H^{-2}(\Omega)), \\ \\
\log(\vt) \in L^\infty(0,T; L^1(\Omega)) \end{array} \right\},
\eF
and
\bFormula{m2}
\Grad \vu \in L^2((0,T) \times \Omega;  R^{3 \times 3}), \
\tn{Q} \in L^2(0,T; W^{2,2}(\Omega; \TS)), \ \Grad \vt \in L^q((0,T) \times \Omega; R^3)
\eF
for any $q < 5/3$.

%%%%%%%%%%%%%%%%%%%%%%%%%%%%%%%%%%%%%%%%%%%%%%%%%%%%%%%%%%%%%%%%%%%%%%%%%%%%%%%%%%%%%%%%%%%%%%%%%%%%%%%%%%

\subsection{Weak solutions}
\label{ws}

The weak solutions are defined in the standard way. Given the anticipated regularity of the velocity field, the
incompressibility constraint (\ref{i21}) makes sense a.e. in the set $(0,T) \times \Omega$, while the momentum balance
(\ref{i22}) is replaced by a family of integral identities
\bFormula{m3}
\int_0^T \intO{ \Big[ \vu \cdot \partial_t \varphi + (\vu \otimes \vu) : \Grad \varphi \Big] } \ \dt =
\int_0^T \intO{ \left( \sigma : \Grad \varphi - \vc{g} \cdot \varphi \right) } \ \dt - \intO{
\vu_0 \cdot \varphi(0, \cdot) }
\eF
satisfied for any test function $\varphi \in \DC([0,T) \times \Omega; R^3)$ and where $\sigma$ is defined in \eqref{i25}.
Here and in what follows, we always
tacitly assume that all quantities appearing under the integrals are (at least)
summable in $(0,T) \times \Omega$.

Similarly, the evolutionary equation (\ref{i23}) for the $Q$-tensor is replaced by
\bFormula{m4}
\int_0^T \intO{ \Big[ \tn{Q} : \partial_t \varphi + [\vu \tn{Q}] : \Grad \varphi + \tn{S}(\Grad \vu, \tn{Q}) : \varphi \Big] }
\ \dt = - \int_0^T \intO{ \Gamma(\vt) \tn{H} : \varphi } \ \dt - \intO{ \tn{Q}_0 : \varphi (0, \cdot) }
\eF
for any $\varphi \in \DC([0,T) \times \Omega; R^{3 \times 3})$, where $\tn{S}(\Grad \vu, \tn{Q})$ is defined in \eqref{i5} and
$\tn{H}$ in \eqref{i20}.

The total energy balance (\ref{i24}) is satisfied in the sense of integral identity
\bFormula{m5}
\int_0^T \intO{ \left[ \left( \frac{1}{2}|\vu|^2 + e \right) \partial_t \varphi + \left( \frac{1}{2}|\vu|^2 + e \right)
\vu \cdot \Grad \varphi + \vc{q} \cdot \Grad \varphi \right] } \ \dt
\eF
\[
= \int_0^T \intO{ \Big[ \sigma \vu \cdot \Grad \varphi + \Gamma(\vt) (\Grad  \tn{Q} : \tn{H}) \cdot \Grad \varphi
- \vc{g} \cdot \vu \varphi \Big] } \ \dt - \intO{ \left( \frac{1}{2}|\vu_0|^2 + e_0 \right) \varphi(0, \cdot) }
\]
for any $\varphi \in \DC([0,T) \times \Omega)$, where $\vc{q}$ is defined in \eqref{i26} and $e$ in \eqref{i27} and we have set
\[
e_0 = \frac{1}{2} |\Grad \tn{Q}_0|^2 + f(\tn{Q}_0)-
\Big( U(\vt_0) - \vt_0 U'(\vt_0) \Big)
G(\tn{Q}_0) + \vt_0.
\]

%%%%%%%%%%%%%%%%%%%%%%%%%%%%%%%%%%%%%%%%%%%%%%%%%%%%%%%%%%%%%%%%%%%%%%%%%%%%%%%%%%%%%%%%%%%%%%%%%%%%%%%%%%5

\subsection{Main result}

We are ready to state the main result of this paper.

\greybox{

\bTheorem{m1}
Let the initial data $\vu_0$, $\tn{Q}_0$, and $\vt_0$ be given such that
\bFormula{m6}
\vu_0 \in L^2(\Omega; R^3), \ \Div \vu_0 = 0, \left\{ \begin{array}{c} \tn{Q}_0 \in W^{1,2}(\Omega; \TS),\\ \\
f(\tn{Q}_0) \in L^1(\Omega) \end{array} \right\}, \ \vt_0 \in L^\infty(\Omega),\
{\rm ess} \inf_\Omega \vt_0 = \underline{\vt} > 0.
\eF
Suppose that the functions $U$, $G$ satisfy the hypotheses (\ref{i3}--\ref{i3bis}),
(\ref{i4}), the transport coefficients
$\mu$, $\kappa$, and $\Gamma$ comply with (\ref{i7}), and that
\[
\vc{g} \in L^\infty (0,T; L^2(\Omega; R^3)).
\]

Then the problem (\ref{i21}--\ref{i24}), (\ref{i29}), (\ref{i30}) admits a weak solution $\vu$, $\tn{Q}$,
$\vt$ in $(0,T) \times \Omega$ in the sense specified in Section \ref{ws}. In addition, there exist positive constants
$c$ and $\lambda$ such that
\bFormula{m7}
\vt(t, \cdot) \geq c \exp (- \lambda t) \underline{\vt} \ \mbox{for all}\ t > 0,
\eF
and the entropy inequality (\ref{i28}) holds in the sense of distributions.
\eT
}

The rest of the paper is devoted to the proof of Theorem \ref{Tm1}. For the sake of simplicity, we set
$\vc{g} = 0$ as the proof in the more general case requires only straightforward modifications.

%%%%%%%%%%%%%%%%%%%%%%%%%%%%%%%%%%%%%%%%%%%%%%%%%%%%%%%%%%%%%%%%%%%%%%%%%%%%%%%%%%%%%%%%%%%%%%%%%%%%%%%%%%5

\section{A priori bounds}
\label{a}

{\it A priori}\/ bounds are natural (formal) constraints imposed on hypothetical smooth solutions
by the equations and the initial data.

%%%%%%%%%%%%%%%%%%%%%%%%%%%%%%%%%%%%%%%%%%%%%%%%%%%%%%%%%%%%%%%%%%%%%%%%%%%%%%%%%%%%%%%%%%%%%%%%%%%%%%%%%%5

\subsection{Energy bounds}

Uniform boundedness in time of the total energy is straightforward consequence of (\ref{i24}). In accordance with
(\ref{m1}), we get
\bFormula{a1}
\vu \in L^\infty(0,T; L^2(\Omega; R^3)),
\eF
\bFormula{a2}
f(\tn{Q}) \in L^\infty (0,T;L^1(\Omega)), \ \mbox{in particular}, \
\tn{Q} \in L^\infty ((0,T) \times \Omega; \TS), \ \Grad \tn{Q} \in L^\infty(0,T; L^2(\Omega; R^{27})),
\eF
and
\bFormula{a3}
\vt \in L^\infty(0,T; L^1(\Omega)),
\eF
where we have anticipated the fact that the absolute temperature is a positive quantity.

%%%%%%%%%%%%%%%%%%%%%%%%%%%%%%%%%%%%%%%%%%%%%%%%%%%%%%%%%%%%%%%%%%%%%%%%%%%%%%%%%%%%%%%%%%%%%%%%%%%

\subsection{Entropy bounds}

Integrating the entropy inequality (\ref{i28}) and using (\ref{a3}), (\ref{i7})$(b)$ we infer that
\bFormula{a4}
\log (\vt) \in L^\infty(0,T; L^1(\Omega)) \cap L^2(0,T; W^{1,2}(\Omega)).
\eF

%%%%%%%%%%%%%%%%%%%%%%%%%%%%%%%%%%%%%%%%%%%%%%%%%%%%%%%%%%%%%%%%%%%%%%%%%%%%%%%%%%%%%%%%%%%%%%%%%%%

\subsection{Bounds based on energy dissipation}
\label{sub:bounds}

Multiplying the entropy inequality (\ref{i28}) by $\vt$ we deduce the \emph{thermal energy balance} in the form
\[
\partial_t \vt + \vu \cdot \Grad \vt - \Div \Big( \kappa(\vt) \Grad \vt \Big) \geq
- \vt \left[ \partial_t \Big( U'(\vt) G(\tn{Q})\Big) + \vu \cdot \Grad \Big( U'(\vt) G(\tn{Q})\Big) \right]
\]
\[
+ \frac{\mu(\vt)}{2} \left| \Grad \vu + \Grad^t \vu \right|^2 + \Gamma(\vt) \tn{H} : \tn{H},
\]
where, furthermore,
\[
\vt  \partial_t \Big( U'(\vt) G(\tn{Q})\Big) = \partial_t \left[ \Big( \vt U'(\vt) - U(\vt) + U(0) \Big) G(\tn{Q}) \right]
+ \Big( U(\vt) - U(0) \Big) \partial_t G(\tn{Q}),
\]
whence
\bFormula{a5}
\partial_t \left[ \Big( \vt U'(\vt) - U(\vt) + U(0) \Big) G(\tn{Q}) + \vt \right] + \vu \cdot \Grad
\left[ \Big( \vt U'(\vt) - U(\vt) + U(0) \Big) G(\tn{Q}) + \vt \right]
\eF
\[
 - \Div \Big( \kappa(\vt) \Grad \vt \Big) \geq \Big( U(0) - U(\vt) \Big) \Big[ \partial_t G(\tn{Q}) +
 \vu \cdot \Grad G(\tn{Q}) \Big] + \frac{\mu(\vt)}{2} \left| \Grad \vu + \Grad^t \vu \right|^2 + \Gamma(\vt) \tn{H} : \tn{H}
\]
\[
= \Big( U(0) - U(\vt) \Big) \mathcal{L} \left[ \frac{\partial G(\tn{Q})}{\partial \tn{Q}} \right]
: \Big[ \tn{S}(\Grad \vu, \tn{Q} ) + \Gamma(\vt) \tn{H} \Big] + \frac{\mu(\vt)}{2} \left| \Grad \vu + \Grad^t \vu \right|^2
+ \Gamma(\vt) \tn{H} : \tn{H}.
\]

Since we already know that $\tn{Q}$ is uniformly bounded, we deduce from (\ref{a5}) that
\bFormula{a6}
\partial_t \intO{ \left[ \Big( \vt U'(\vt) - U(\vt) + U(0) \Big) G(\tn{Q}) + \vt \right] }
\eF
\[
\geq \intO{ \left[ \frac{\mu(\vt)}{4} \left| \Grad \vu + \Grad^t \vu \right|^2 + \frac{\Gamma(\vt)}{2} \tn{H} : \tn{H} -
C\Big| U(0) - U(\vt) \Big|^2 \right]},
\] where $C>0$ is an explicitly computable constant depending only on $\|\tn{Q}\|_{L^\infty}$, $G$ and $\xi$.

By virtue of hypotheses (\ref{i3}--\ref{i3bis}), the function $U'(\vt)$ is bounded,
therefore we may combine (\ref{a6})
with the energy estimate (\ref{a3}), use assumptions (\ref{i7}) and the convexity of $f$ to conclude that
\bFormula{a7}
\Grad \vu \in L^2((0,T) \times \Omega; R^3), \ \tn{Q} \in L^2(0,T; W^{2,2}(\Omega; \TS)), \ \mbox{and}\
\mathcal{L} \left[ \frac{\partial f(\tn{Q})}{\partial \tn{Q}} \right] \in  L^2(0,T; L^2(\Omega; \TS)).
\eF
Finally, thanks to (\ref{a1}--\ref{a2}) and \eqref{a7},
a comparison of terms in \eqref{i23} gives also
\bFormula{a7x}
  \partial_t \QQ \in L^1(0,T;L^3(\Omega; \TS)).
\eF
%

%%%%%%%%%%%%%%%%%%%%%%%%%%%%%%%%%%%%%%%%%%%%%%%%%%%%%%%%%%%%%%%%%%%%%%%%%%%%%%%%%%%%%%%%%%%%%%%%%%%%%%%%%%%%%%%%%%%%%5

\subsection{Bounds on the temperature gradient}

The estimates on the temperature gradient are obtained by multiplying (\ref{a5}) on
$- ( 1 + \vt)^{-\alpha}$, $\alpha > 0$. Note that
\[
  \frac{1}{( 1 + \vt)^\alpha } \Div \Big( \kappa(\vt) \Grad \vt \Big)
   = \Div \left( \frac{\kappa(\vt)}{(1 + \vt)^\alpha } \Grad \vt \right)
    + \frac{4 \alpha}{ (1-\alpha)^2 } \kappa(\vt)
        \left| \Grad (1 + \vt)^{\frac{1 - \alpha}{2} } \right|^2.
\]
Observe also that the right hand side of \eqref{a5} is uniformly bounded in $L^1$
thanks to the previous estimates. Moreover, as we integrate by parts the terms depending
on $U$, we obtain, on the right hand side, the terms
\[
  G(\QQ) \frac{\vt}{(1 + \vt)^\alpha}U''(\vt) \vt_t
   + \frac{\vt}{(1 + \vt)^\alpha} U'(\vt) G(\QQ)_t,
\]
which we need to control (as well as similar quantities
depending on the transport part of the material derivative, which can
be treated in the same way). Integrating by parts in time, we get
\bFormula{gh1}
  \partial_t \big( G(\QQ) Y_\alpha(\vt) \big)
   + \left[ \frac{\vt}{(1 + \vt)^\alpha} U'(\vt) - Y_\alpha(\vt) \right] G(\QQ)_t,
\eF
where
\[
  Y_\alpha(\vt) \equiv \int_1^{\vt} \frac{s}{(1 + s)^\alpha}U''(s) \ \ds
\]
and, thanks to (\ref{i3}--\ref{i3bis}), the function in square
brackets in \eqref{gh1} goes like $\vt^{1/2 - \alpha}$ for
large $\vt$. Hence, due to \eqref{a3}, it lies (at least) in
$L^\infty(0,T;L^2(\Omega))$.

Now, using that $G$ is bounded with its first derivatives and
recalling \eqref{a7x}, all terms in \eqref{gh1} can be controlled. Hence,
after a straightforward manipulation we can
conclude that
\bFormula{a8}
\Grad (1 + \vt)^{\frac{1 - \alpha}{2} } \in L^2((0,T) \times \Omega;R^3) \ \mbox{for any}\ \alpha > 0.
\eF

The a-priori bounds obtained here are enough to make the weak formulation (\ref{m3}--\ref{m5}) meaningful. In particular,
 the pressure $p$ (cf. \eqref{i25}) can be ``computed'' directly from \eqref{m3}
 (cf.~\cite{MAL2} and \cite{FFRS} for more details) and it is possible to
 obtain
 \[
 p\in L^{5/3}((0,T)\times \Omega)\,.
 \]

It can be shown that the {\it a priori} bounds obtained in this section are strong enough in order to establish the weak sequential
stability of the family of solutions to our problem. However, we do not pursue this path and pass directly to the construction of a
family of \emph{approximate} solutions.

%%%%%%%%%%%%%%%%%%%%%%%%%%%%%%%%%%%%%%%%%%%%%%%%%%%%%%%%%%%%%%%%%%%%%%%%%%%%%%%%%%%%%%%%%%%%%%%%%%%%%%%%%%%%%%%%%%%%%%%

\section{Approximate problems}
\label{p}

The weak solution, the existence of which is claimed in Theorem \ref{Tm1},
will be constructed by means of a family of approximate problems.

%%%%%%%%%%%%%%%%%%%%%%%%%%%%%%%%%%%%%%%%%%%%%%%%%%%%%%%%%%%%%%%%%%%%%%%%%%%%%%%%%%%%%%%%%%%%%%%%%%%%%

\subsection{Approximate velocity fields}
\label{pp}

The velocity field $\vu$ is obtained via the standard
Faedo-Galerkin method based on the family of finite dimensional spaces
\[
X_N = \left\{ \vc{v} \in C^\infty(\Omega; R^3) \ \Big| \ \vc{v}
 - \ \mbox{a trigonometric polynomial of order}\ N, \ \Div \vc{v} = 0 \right\}.
\]

Accordingly, the momentum equation (\ref{m3}) is replaced by a finite system of ordinary differential equations
\bFormula{p1}
\frac{{\rm d}}{{\rm d}t} \intO{ \vu \cdot \vc{v} }
\eF
\[
= \intO{ [\vu]_\delta \otimes \vu : \Grad \vc{v} } - \delta \intO{
|\Grad \vu |^{r-2} \Grad \vu : \Grad \vc{v} } - \intO{ \mu(\vt) (
\Grad \vu + \Grad^t \vu) : \Grad \vc{v} }
\]
\[
+ \intO{ \Big( \Grad \tn{Q} \odot \Grad \tn{Q} \Big) : \Grad \vc{v}
}
\]
\[
- \intO{ \left\{  2 \xi \left( \tn{H}_{m,\delta} : \tn{Q} \right) \left(
\tn{Q} + \frac{1}{3} \tn{I} \right) - \xi \left[ \tn{H}_{m,\delta} \left(
\tn{Q}  + \frac{1}{3} \tn{I} \right) + \left( \tn{Q}  + \frac{1}{3}
\tn{I} \right) \tn{H}_{m,\delta} \right] + \left( \tn{Q} \tn{H}_{m} - \tn{H}_{m,\delta}
\tn{Q} \right) \right\} : \Grad \vc{v} },
\]
with the initial condition
\bFormula{p2}
\intO{ \vu(0, \cdot) \cdot \vc{v} } =  \intO{ [\vu_0]_\delta \cdot \vc{v} }
\eF
for any $\vc{v} \in X_N$,
where
\bFormula{p3}
\tn{H}_{m,\delta} = \Delta \tn{Q} - \mathcal{L} \left [
\frac{\partial f_m (\tn{Q})} {\partial \tn{Q}} \right] + U_\delta(\vt) \mathcal{L} \left [
\frac{\partial G (\tn{Q})} {\partial \tn{Q}} \right].
\eF

Here $r\in (3,10/3)$, $\delta$ and $m$ stand for positive parameters; $[ \vu ]_{\delta}$ denotes the standard
regularization with respect to the $x$-variable by means
of a family of convolutions, while $\{ f_m \}_{m > 0}$ is a family of smooth convex functions defined on $\TS$ such that
\[
f_m \leq f\ \hbox{ for all } m\in\mathbb{N},\quad
f_{m_1}(\tn{Q})\le f_{m_2}(\tn{Q}),\hbox{ for all } m_1\le
m_2, \hbox{ for all }\tn{Q}\in\TS,
\]
\[
\left\{ \begin{array}{c} f_m \to f,\ \mbox{uniformly on compact subsets of}\ \mathcal{D}[f], \\ \\
\mathcal{L} \left[ \frac{ \partial f_m }{\partial \QQ } \right] \to
\mathcal{L} \left[ \frac{ \partial f }{\partial \QQ } \right],\ \mbox{uniformly on compact subsets of}\ \mathcal{D}[f],\\ \\
f_m \to \infty \ \mbox{uniformly in} \ \TS \setminus \mathcal{D}[f] \end{array} \right\}\ \mbox{as}\ m \to \infty,
\]
\[
c^1_m |\tn{Q}| - c^2_m \leq \left| \mathcal{L} \left[ \frac{ \partial f_m }{\partial \tn{Q} } \right] \right| \leq
C^1_m |\tn{Q}| + C^2_m \ \mbox{for all} \ \tn{Q} \in \TS, \ m > 0.
\]
Notice that the term $\delta |\Grad\vu|^{r-2}\Grad\vu$ guarantees additional
regularity of $\vu$ needed in the $Q$-tensor equation (cf.~estimate \eqref{p13}
and Section~\ref{sec:N} below).

%Finally, $\tilde G$ denotes a truncation of $G$, specifically,
%
%\bFormula{Gtrunc}
 % \tilde G \in C^1_c (\TS), \ \tilde G(\tn{Q}) \geq 0, \ \
  %\tilde G (\tn{Q}) = G(\tn{Q}) \ \mbox{for all}\ \tn{Q} \in \mathcal{D}[f], \,\exists\, M>0\,:\, \tilde G \hbox{ is constant on }\{|\tn{Q}|\geq M\} .
%\eF

%%%%%%%%%%%%%%%%%%%%%%%%%%%%%%%%%%%%%%%%%%%%%%%%%%%%%%%%%%%%%%%%%%%%%%%%%%%%%%%%%%%%%%%%%%%%%%%%%%%%%

\subsection{$Q$-tensors}

The equation governing the time evolution of the approximate $Q$-tensors
reads
\bFormula{p4}
\partial_t \tn{Q} + (\vu \cdot \Grad) \tn{Q} - \tn{S} (\Grad \vu, \tn{Q} ) =[\Gamma (\vt )]_\ep
\Big( \Delta \tn{Q} - \mathcal{L} \left [
\frac{\partial f_m (\tn{Q})} {\partial \tn{Q}} \right] + U_\delta (\vt) \mathcal{L} \left [
\frac{\partial G (\tn{Q})} {\partial \tn{Q}} \right]
\Big),
\eF
with
\bFormula{p5}
\tn{Q}(0, \cdot) = [\tn{Q}_0]_{\delta}.
\eF
Here, $[\Gamma (\vt) ]_{\ep}$ denotes a regularization (via convolutions)
of $\Gamma(\vt)$ with respect to both $t$ and $x$ variables.
Furthermore, $[\tn{Q}_0]_\delta \in C^\infty(\Omega; \TS)$, such that
\bFormula{p5b}
  [\tn{Q}_0]_\delta \to \tn{Q}_0 \ \mbox{a.a. in} \ \Omega,\
  f([\tn{Q}_0]_\delta ) \to f(\tn{Q}_0) \ \mbox{in} \ L^1(\Omega).
\eF
An explicit construction of $[\tn{Q}_0]_\delta$ can be obtained, for
instance, by truncation and mollification (recall that the domain of $f$
is an open set). Then, the convergence property in \eqref{p5b}
can be verified using the dominated convergence theorem.

Finally, $U_\delta\,:\,\RR\to \RR$ is a bounded truncation of $U$ satisfying (\ref{i3}--\ref{i3bis}) and such that
$U_\delta'(\vt)=U_\delta'(0)$ for $\vt\leq 0$.

%%%%%%%%%%%%%%%%%%%%%%%%%%%%%%%%%%%%%%%%%%%%%%%%%%%%%%%%%%%%%%%%%%%%%%%%%%%%%%%%%%%%%%%%%%%%%%%%%%%%%

\subsection{Thermal energy balance}

The approximate temperature is determined via a ``heat'' equation of the form
\bFormula{p6}
\partial_t \vt  + \vu \cdot \Grad \vt  - \Div \Big( \kappa(\vt) \Grad \vt \Big)
\eF
\[
= - \vt \partial_t \Big( U'_\delta (\vt)  G (\tn{Q}) \Big) - \vt \vu
\cdot \Grad \Big( U'_\delta(\vt) G (\tn{Q} ) \Big) +
\frac{\mu(\vt)}{2} \Big| \Grad \vu + \Grad^t \vu \Big|^2
\]
\[
+ [\Gamma(\vt)]_\ep \tn{H}_{m,\delta}: \tn{H}_{m,\delta}
+ \delta |\Grad \vu |^r,
\]
where $\tn{H}_{m,\delta}$ is defined in \eqref{p3}.
%
%\[
%\tn{H}_{m,\delta}:=
%\left( \Delta \tn{Q} - \mathcal{L}
%\left [ \frac{\partial f_m (\tn{Q})} {\partial \tn{Q}} \right]
%+ U_\delta(\vt) \mathcal{L} \left [ \frac{\partial \tilde G (\tn{Q})} {\partial \tn{Q}} \right]
%\right).
%\]
%
Equation \eqref{p6} is complemented with the initial condition
\bFormula{p7}
  \vt(0, \cdot) = [\vt_0]_{\delta},
\eF
where, similarly to (\ref{p2}), $[\vt_0]_\delta$ denotes a regularization in the space variables.

%%%%%%%%%%%%%%%%%%%%%%%%%%%%%%%%%%%%%%%%%%%%%%%%%%%%%%%%%%%%%%%%%%%%%%%%%%%%%%%%%%%%%%%%%%%%%%%%%%%%%

\subsection{Existence of approximate solutions and uniform bounds}

Our program for the remaining part of the paper will be to construct approximate solutions to the problem
(\ref{p1}--\ref{p7}) and let
successively
\[
  m \to \infty, \ N \to \infty, \ \ep \to 0, \ \mbox{and, finally,} \ \delta \to 0
\]
in order to recover in the limit a weak solution to the problem (\ref{i21}--\ref{i24}),
(\ref{i29}--\ref{i30}), the existence of which is claimed in Theorem \ref{Tm1}.

For fixed values of the parameters $m$, $N$, $\ep$, and $\delta$, we can construct \emph{local-in-time}
solutions to the approximate system by means of a Schauder fixed-point argument. This procedure is
similar to the one sketched in \cite[Section~5]{FFRS} (see also \cite[Chapter~6]{FEINOV} for further details);
hence we leave it to the reader. Moreover, the local solutions can be extended to the whole time interval
$[0,T]$ as soon as suitable uniform estimates analogous to the {\it a priori} bounds are established.
Actually, to simplify notations, we shall directly assume that solutions are defined on the
whole $(0,T)$ already in the approximation.

%%%%%%%%%%%%%%%%%%%%%%%%%%%%%%%%%%%%%%%%%%%%%%%%%%%%%%%%%%%%%%%%%%%%%%%%%%%%%%%%%%%%%%%%%%%%%%%%%%%%%

\subsubsection{Energy bounds}

Bounds on the total energy are obtained in the same way as the {\it a priori} bounds. We take $\vc{v} = \vu(t, \cdot)$
 as a test function in (\ref{p1}), multiply (\ref{p4}) by $- \tn{H}_{m, \delta}$, and add the resulting expression to (\ref{p6}) to obtain
\bFormula{p8}
\frac{{\rm d}}{{\rm d}t} \intO{ \left[ \frac{1}{2} |\vu|^2 +
\frac{1}{2} |\Grad \tn{Q} |^2 + f_m (\tn{Q}) - \Big( U_\delta(\vt) - \vt U'_\delta(\vt) \Big)
  G (\tn{Q}) + \vt \right] } = 0.
\eF
We deduce that
\bFormula{p9}
{\rm ess} \sup_{t \in (0,T)} \| \vu(t, \cdot) \|_{L^2(\Omega; R^3)} \leq c,
\eF
\bFormula{p10}
{\rm ess} \sup_{t \in (0,T)} \| \tn{Q} (t, \cdot) \|_{W^{1,2}(\Omega; \TS)} \leq c,
\eF
and
\bFormula{p11}
{\rm ess} \sup_{t \in (0,T)} \| \vt (t, \cdot) \|_{L^1(\Omega)} \leq c.
\eF

In addition, since $\vu$ ranges in the finite dimensional space $X_N$ consisting of smooth functions, we get
\bFormula{p12}
\sup_{t \in [0,T]} \| \vu(t, \cdot) \|_{C^k (\Omega; R^3)} \leq c(k,N) \ \mbox{for any} \ k = 0,1,\dots
\eF

%%%%%%%%%%%%%%%%%%%%%%%%%%%%%%%%%%%%%%%%%%%%%%%%%%%%%%%%%%%%%%%%%%%%%%%%%%%%%%%%%%%%%%%%%%%%%%%%%%%%%%%%%%%

\subsubsection{Bounds on the $Q$-tensors}

Since the approximate $Q$-tensors satisfy equation (\ref{p4}), where the leading coefficient $[\Gamma(\vt)]_\ep$
is smooth and the coupling term $\vu \cdot \Grad \tn{Q} - \tn{S} (\Grad \vu, \tn{Q})$ is regular thanks
to \eqref{p12}, we may bootstrap the maximal regularity
estimates of $L^q$-type (see e.g. Krylov \cite{Kry})
to deduce that
\bFormula{p13}
\| \partial_t \tn{Q} \|_{L^q(0,T; L^q(\Omega; \TS))} + \| \tn{Q} \|_{L^q(0,T; W^{2,q}(\Omega; \TS))} \leq c(q,N,m,\ep, \delta)
\ \mbox{for any} \ 1 \leq q < \infty.
\eF

Now, we can go back to (\ref{p1}) and to use (\ref{p13}) to conclude that
\bFormula{p14}
{\rm ess} \sup_{t \in (0,T)} \| \partial_t \vu(t, \cdot) \|_{C^k(\Omega;R^3)} \leq c(k,N,m,\ep, \delta).
\eF

%%%%%%%%%%%%%%%%%%%%%%%%%%%%%%%%%%%%%%%%%%%%%%%%%%%%%%%%%%%%%%%%%%%%%%%%%%%%%%%%%%%%%%%%%%%%%%%%%%%%%%%%%%%

\subsubsection{Strict positivity of the absolute temperature}

This is one of the most delicate steps in the proof of Theorem \ref{Tm1}. Our aim is to prove
that $\vt$ is strictly positive already at the approximate level and {\it uniformly}\/
with respect to all approximation parameters. To this aim, we apply the parabolic comparison
theorem to equation (\ref{p6}) written in the form
\bFormula{p15}
\Big( 1 + \vt U''_{\delta}(\vt)  G(\tn{Q}) \Big) \Big( \partial_t \vt  + \vu \cdot \Grad \vt \Big)
 - \Div \Big( \kappa(\vt) \Grad \vt \Big)
\eF
\[
 = - \vt  U'_\delta (\vt) \mathcal{L} \left[ \frac{ \partial  G (\tn{Q})}{\partial \tn{Q}} \right]
  : \Big( \tn{S}(\Grad \vu, \tn{Q}) + [\Gamma(\vt)]_\ep \tn{H}_{m,\delta} \Big)
  + \frac{\mu(\vt)}{2} \Big| \Grad \vu + \Grad^t \vu \Big|^2
\]
\[
+ [\Gamma(\vt)]_\ep \tn{H}_{m,\delta} : \tn{H}_{m, \delta}
+ \delta |\Grad \vu |^r.
\]

Now, a short inspection of (\ref{i5}) yields that all terms can be expressed as products of
scalar quantities depending on $\tn{Q}$ and the symmetric gradient $\ep(\vu)$ with the only exception of the commutator
\[
\omega(\vu) \tn{Q} - \tn{Q} \omega(\vu).
\]
Fortunately, thanks to assumption \eqref{i4}, we have
\[
\mathcal{L} \left[ \frac{ \partial  G (\tn{Q})}{\partial \tn{Q}} \right]: \Big[ \omega(\vu) \tn{Q} - \tn{Q} \omega(\vu) \Big] =
2 \left( \mathcal{L} \left[ \frac{ \partial  G (\tn{Q})}{\partial \tn{Q}} \right] \tn{Q} \right) : \omega(\vu) = 0.
\]
Let us now observe that, by (\ref{i3}--\ref{i3bis}), there exists $c>0$ such that $|U'(\vt)\vt|\le c \vt^{1/2}$
for all $\vt\ge 0$, and the same inequality can be required to hold for the approximations $U_\delta$.
Moreover, we can suppose $G$ bounded because the domain of $f$ is bounded and 
consequently the behavior of $G$ for large $|\tn{Q}|$ is not relevant. On account 
of these considerations, it is not difficult to arrive at
\bFormula{p16}
  \Big( 1 + \vt U''_{\delta}(\vt)  G(\tn{Q}) \Big) \Big( \partial_t \vt  + \vu \cdot \Grad \vt \Big)
   - \Div \Big( \kappa(\vt) \Grad \vt \Big)
\eF
\[
\geq - \Lambda \vt +
\frac{\mu(\vt)}{4} \Big| \Grad \vu + \Grad^t \vu \Big|^2
+ \frac{[\Gamma(\vt)]_\ep}{2} \tn{H}_{m,\delta} : \tn{H}_{m,\delta}
+ \delta |\Grad \vu |^r,
\]
where $\Lambda$ is a positive constant depending only on the structural
properties of the functions $U_\delta$ and $ G$.

Testing \eqref{p15} by $\vt^-$ and using the fact that $U_\delta''=0$ for $\vt<0$, we obtain, by standard arguments, the
non-negativity of $\vt$.

Next, multiplying equation (\ref{p15}) by  $\kappa(\vt) \partial_t \vt$ and making use of the available bounds (\ref{p11}), (\ref{p12}), (\ref{p13}),
and \eqref{i3bis}, we deduce
\bFormula{p17}
  {\rm ess} \sup_{t \in (0,T)} \| \vt \|_{W^{1,2}(\Omega)} \leq c(N,m,\ep, \delta), \
    \| \partial_t \vt \|_{L^2((0,T) \times \Omega)} \leq c(N,m,\ep, \delta).
\eF
Thus, introducing
\[
  \Theta = \mathcal{K}(\vt), \ \mathcal{K}' = \kappa, \ \mathcal{K}(0)=0,
\]
and noting that the coefficient $(1 + \vt U''_{\delta}(\vt)  G(\tn{Q}))$ is uniformly
bounded, applying standard elliptic regularity results to~\eqref{p15} we arrive at
\bFormula{p18new}
  \| \Theta \|_{L^2(0,T; W^{2,2}(\Omega))} \leq c(N,m,\ep, \delta).
\eF
whence, by \eqref{p17} and interpolation, it is not difficult to
obtain
\bFormula{p18new2}
  \| \vt \|_{L^2(0,T; W^{2,3/2}(\Omega))} \leq c(N,m,\ep, \delta).
\eF
%
%\bFormula{p18}
%\| \vt \|_{L^2(0,T; W^{2,2}(\Omega))} \leq c(N,m,\ep, \delta).
%\eF
%
Next, we may rewrite (\ref{p16}) in the form
\bFormula{p19}
a(t,x) \Big( \partial_t \Theta + \lambda \Theta +\vu \cdot \Grad \Theta \Big) - \Delta \Theta
 = g \geq 0,
\eF
where the previous estimates imply that $g\in L^p((0,T) \times \Omega)$ for all $p\in[1,\infty)$
and we have set
\[
a = \frac{ 1 + \vt U''_{\delta}(\vt)  G(\tn{Q}) }{\kappa(\vt)} \in L^\infty((0,T) \times \Omega),
 \ a\geq \underline{a}>0,
\]
where $\lambda > 0$ depends only on $\Lambda$ and the structural properties of the
function $\kappa$, and $\underline{a}$ depends only
on $\Ov\kappa$.
For instance, one may take
\bFormula{lamb}
  \lambda \equiv {\rm ess} \inf \frac{\Lambda \vt}{\mathcal{K}(\vt) a} > 0.
\eF
Then, multiplying \eqref{p19} by $e^{\lambda t}$ and setting
$z\equiv e^{\lambda t} \Theta$, we get
\bFormula{p19b}
  a(t,x) \Big( z_t + \vu \cdot \Grad z \Big) - \Delta z
 = g e^{\lambda t} \geq 0,
\eF
whence, testing by $ - ( z - \underline{\Theta})^{-} $, 
where ${\underline\Theta} = \mathcal{K} ({\underline\vt})$,
and integrating by parts the terms with $a$, we obtain the desired conclusion
\bFormula{p20}
  \vt (t, \cdot) \geq c \exp(-\lambda t) \underline{\vt} \ \mbox{for all} \ t \in [0,T].
\eF
Let us note that in principle estimate (\ref{p20}) holds for smooth $a$ (see e.g. \cite{FRIE1}).
However, considering the class of solutions specified through (\ref{p17}), (\ref{p18new2}),
the procedure can easily be extended to any bounded measurable $a$ by
means of density arguments.
Indeed, we can consider a smooth (H\"older continuous) approximation $a_h$ of $a$
and a smooth approximation $g_h$ of $g$,
and we can assume $a_h\to a$, $g_h\to g$ in $L^p$ for every $p\in (1,\infty)$.
Then, by the standard maximum principle argument (cf. \cite{FRIE1}) we can conclude that there exists
$\vt_h\in H^1(0,T;L^2(\Omega))\cap L^2(0,T; W^{2,2}(\Omega))$
(in fact, $\vt_h$ will be even smoother) satisfying \eqref{p20}. Then, also the limit $\widetilde{\vt}$
of $\vt_h$ satisfies the same inequality. This concludes the proof
because we have uniqueness of solutions for \eqref{p19},
which can be obtained simply testing the difference of the two
equations by the time derivative of the difference of two solutions.
An analogous proof of the maximum principle argument for \eqref{p19}
with bounded measurable $a$ can be found, e.g., in  \cite[Prop.~3.6, p.~293]{KRROSP}.

Estimate (\ref{p20}) coincides with (\ref{m7}) claimed in Theorem \ref{Tm1}. It is
remarkable that (\ref{p20}) is independent of the parameters
$m$, $N$, $\ep$, and $\delta$. Indeed, the choice \eqref{lamb}
of $\lambda$ is independent of all approximations.

%%%%%%%%%%%%%%%%%%%%%%%%%%%%%%%%%%%%%%%%%%%%%%%%%%%%%%%%%%%%%%%%%%%%%%%%%%%%%%%%%%%%%%%%%%%%%%%%%%%%%%%%%%%%

\subsubsection{Moser and regularity estimates on the absolute temperature}

We prove here that the absolute temperature is (globally in time) H\"older continuous
at least at the approximate level. Although this fact is essentially a consequence of
well-known techniques for parabolic equations (see, e.g., \cite{LADSOLUR}), in view of the
fact that \eqref{p15} depends on $\vt$ in a somehow intricated way we give, for
the convenience of the reader, at least the highlights of a direct proof
based on Moser iterations.

First of all, we rewrite \eqref{p15} in the form
\bFormula{p15b}
  \Big( 1 + \vt U''_{\delta}(\vt)  G(\tn{Q}) \Big) \Big( \partial_t \vt  + \vu \cdot \Grad \vt \Big)
  - \Div \Big( \kappa(\vt) \Grad \vt \Big)
   = - \vt  U'_\delta (\vt) \ell + \nu,
\eF
where the functions $\ell$ and $\nu$ collect the various quantities on the right hand side
of~\eqref{p15} and satisfy
\bFormula{p15c}
  \| \ell \|_{L^q(0,T; L^q(\Omega))}
   + \| \nu \|_{L^q(0,T; L^q(\Omega))} \le c(q,N,m,\ep,\delta)
   \ \mbox{for any} \ 1 \leq q < \infty.
\eF
Then, by the lower bound \eqref{p20} we are allowed to test \eqref{p15c} by
$\vt^{p-1}$ for a generic $p>1$. Integrating with respect to space variables,
we then obtain
\bFormula{p15d}
  \partial_t \intO{ \left[ \left( \frac1p \vt^p + H_{\delta,p}(\vt)  G(\tn{Q}) \right)
    + \frac{4(p-1)}{p^2} \kappa(\vt) \left| \Grad \vt^{p/2} \right|^2 \right] }
\eF
\[
  = \intO{ \left[ H_{\delta,p}(\vt) \mathcal{L} \left[ \frac{\partial  G (\tn{Q}) }{\partial \tn{Q}} \right]
   \big( \partial_t \tn{Q} + \vu \cdot \Grad \tn{Q} \big)
   - \vt^p  U'_\delta (\vt) \ell + \nu \vt^{p-1} \right] },
\]
where we have set
\[
  H_{\delta,p}(\vt):= \int_1^{\vt} r^p U_\delta''(r)\ \dr
\]
and we can notice that, due to (\ref{i3}--\ref{i3bis})
(or, more precisely, the analogue for $U_\delta$), we have
\[
  H_{\delta,p}(\vt) \le \frac{c}{p} \left( 1 + \vt^{p-1/2} \right).
\]
Then, owing to the fact that $H_{\delta,p}(\vt)  G(\tn{Q})$ is nonnegative,
and using the regularity of $\tn{Q}$ and $\vu$ to estimate the right hand side,
we see that \eqref{p15d} assumes in fact the structure
\bFormula{p15e}
  \partial_t \| \vt \|_{L^p(\Omega)}^p
   + c_1 \| \Grad \vt^{p/2} \|_{L^2(\Omega);R^3}^2
   \le c_2 (1 + p) \intO{ \left( \phi \vt^{1/2} \right) \left( 1 + \vt^{p-1} \right) },
\eF
where $c_1,c_2>0$ are absolute constants and
$\phi$ collects all the terms depending on $\vu$ and $\tn{Q}$
(also through $\ell$ and $\nu$) and belongs to $L^q((0,T)\times\Omega)$
for all $q\in[1,\infty)$.
Thus, for the sake of applying Moser iterations, relation
\eqref{p15e} shows that we can proceed by working exactly
as in the linear case (with forcing term given here
by $\phi \vt^{1/2}$). Actually, by \eqref{p17}, \eqref{p18new2}
and interpolation we have $\phi \vt^{1/2} \in L^r((0,T)\times\Omega)$
for
%, say, $r\in [1,20)$ (actually, for Moser iterations it would have been sufficient
%to have $\phi \vt^{1/2} \in L^p((0,T)\times\Omega)$ for
%
a suitable $r>3$. Thus, applying the results in
\cite{LADSOLUR}, we obtain the bound
\bFormula{p15f}
  \| \vt \|_{L^\infty((0,T)\times \Omega)}\le c(N,m,\ep,\delta),
\eF
and this estimate depends on the approximation parameters, but is uniform with respect
to time. Finally, applying the theory developed, e.g., in
\cite[Section~4]{KRYSAF}, we obtain
\bFormula{p15g}
  \| \vt \|_{C^\alpha([0,T]\times \Omega)}\le c(\alpha,N,m,\ep,\delta),
  \ \mbox{for some} \ \alpha>0,
\eF
globally in time.

%%%%%%%%%%%%%%%%%%%%%%%%%%%%%%%%%%%%%%%%%%%%%%%%%%%%%%%%%%%%%%%%%%%%%%%%%%%%%%%%%%%%%%%%%%%%%%%%%%%%%%%%%%%%

\subsubsection{Estimates based on energy dissipation}

The energy bounds analogous to (\ref{a7}), (\ref{a8}) can be deduced from (\ref{p15})
exactly as in Subsection~\ref{sub:bounds}. Actually, although $\tn{Q}$ may not lie in
$L^\infty$ at this stage, the term
\[
  \mathcal{L} \left[ \frac{\partial  G(\tn{Q})}{\partial \tn{Q}} \right]
   :  \tn{S}(\Grad \vu, \tn{Q} )
\]
does not suffer any summability loss since  we can suppose that $ G$ has zero gradient when $| \tn{Q} |$ is large.
Hence, we get
\bFormula{p21}
  \| \Grad \vu \|_{L^2((0,T) \times \Omega; R^{3 \times 3})} \leq c,\
   \delta^{1/r} \| \Grad \vu \|_{L^r((0,T) \times \Omega; R^{3 \times 3})} \leq c,
\eF
\bFormula{p22}
\left\| \mathcal{L} \left[ \frac{\partial f_m (\tn{Q}) }{\partial \tn{Q}} \right] \right\|_{L^2((0,T) \times \Omega; \TS)} \leq c,
\eF
\bFormula{p23}
\| \tn{Q} \|_{L^2(0,T; W^{2,2}(\Omega; \TS))} \leq c,
\eF
%
%  G-N DA SPOSTARE
%
%Hence, applying the standard Gagliardo-Nirenberg inequality,
%\bFormula{p23bis}
%\| \Grad\tn{Q} \|_{L^4((0,T) \times \Omega; R^{27})} \leq c,
%\eF
%
where all constants are independent of the parameters $m$, $N$, $\ep$, and $\delta$.

By \eqref{p10}, \eqref{p23} and interpolation, we also infer
\bFormula{p23b}
  \| \tn{Q} \|_{L^{10}((0,T)\times \Omega; \TS))}
   + \| \Grad \tn{Q} \|_{L^{10/3}((0,T)\times \Omega; R^{27}))} \leq c.
\eF
Analogously, from \eqref{p9} and the first \eqref{p21} we get
\bFormula{p23b2}
  \| \vu \|_{L^{10/3}((0,T)\times \Omega; R^3))} \leq c.
\eF
Coupling \eqref{p23b}, \eqref{p23b2},  and \eqref{p21}, we obtain a uniform bound for the
coupling terms in \eqref{p4}, namely
\bFormula{p23c}
  \big\| \vu \cdot \Grad \tn{Q} - \SSS(\Grad \vu, \tn{Q}) \big\|_{L^{10/7}((0,T)\times \Omega; \TS)} \leq c,
\eF
whence (\ref{p22}--\ref{p23}) and a comparison of terms in \eqref{p4} lead to
\bFormula{p23d}
  \| \partial_t \tn{Q} \|_{L^{10/7}((0,T)\times \Omega; \TS)} \leq c.
\eF
We notice however that estimates (\ref{p23b}--\ref{p23d}) will
be improved in the sequel.

%%%%%%%%%%%%%%%%%%%%%%%%%%%%%%%%%%%%%%%%%%%%%%%%%%%%%%%%%%%%%%%%%%%%%%%%%%%%%%%%%%%%%%%%%%%%%%%%%%%%%%%%%%%%

\section{Convergence to the limit system}
\label{c}

The bounds derived in the previous section, being global in time,
are sufficient for extending the approximate solutions to the desired existence interval.
For this reason, and in order to avoid technicalities, we shall directly
assume that solutions are defined over the whole $(0,T)$ already at the approximate
level. Our ultimate goal is to perform the limits
\[
  m \to \infty, \ N \to \infty, \ \ep \to 0, \ \mbox{and, finally,} \ \delta \to 0.
\]

%%%%%%%%%%%%%%%%%%%%%%%%%%%%%%%%%%%%%%%%%%%%%%%%%%%%%%%%%%%%%%%%%%%%%%%%%%%%%%%%%%%%%%%%%%%%%%%%%%%%%%%%%%%%

\subsection{Uniform bounds on the temperature}
\label{unif:t}

Testing \eqref{p15} by $-(1-\vt)^{-\alpha}$, $\alpha>0$, it is not difficult
to arrive at the analogue of \eqref{a8}, namely
\bFormula{p24}
  \left\| \Grad (1 + \vt)^{\frac{1 - \alpha}{2}}  \right\|_{L^2(0,T; L^2 (\Omega; R^3))} \leq c,
   \ \mbox{for any} \ \alpha > 0,
\eF
Indeed, all terms on the \rhs\ of \eqref{p15} are uniformly bounded in $L^2$. In particular,
$\mathcal{L} \left[ \frac{ \partial G (\QQ)}{\partial \QQ} \right] : \tn{S}(\Grad \vu, \tn{Q})$
has the same regularity of $\Grad \vu$ because we can suppose $G$ to be have zero gradient
in the set where $|\tn{Q}|$ is large. Note also that regularity and strict
positivity of the (approximate) $\vt$ are essential in order that this procedure makes sense.

Coupling the information coming from (\ref{p11}) and (\ref{p24}), which are both independent
of all approximation parameters, and using interpolation arguments together with the Sobolev imbedding
$W^{1,2}(\Omega) \hookrightarrow L^6(\Omega)$, we arrive at
\bFormula{c1}
  \| \vt \|_{L^q((0,T) \times \Omega)} \leq c \ \mbox{for any} \ 1 \leq q < 5/3.
\eF
Moreover, another application of interpolation (see also \cite[Sec.~4]{FFRS})
permits to obtain
\bFormula{c2}
  \| \Grad \vt \|_{L^q((0,T) \times \Omega; R^3)} \leq c \ \mbox{for any}\ 1 \leq q < 5/4.
\eF
Next, thanks to the lower bound established in (\ref{p20}) and to the H\"older regularity
\eqref{p15g}, we can (rigorously) divide the heat equation (\ref{p6}) by $\vt$ to obtain
the entropy relation
\bFormula{c3}
  \partial_t \Big( \log(\vt) + U'_\delta (\vt)  G(\tn{Q}) \Big) + \vu \cdot \Grad \Big( \log(\vt)
   + U'_\delta (\vt)  G(\tn{Q}) \Big) - \Div \left( \frac{\kappa(\vt)}{\vt} \Grad \vt \right)
\eF
\[
  = \frac{1}{\vt} \left( \frac{\mu(\vt)}{2} \left| \Grad \vu + \Grad^t \vu \right|^2
   + \delta |\Grad \vu |^r + [\Gamma(\vt)]_\ep |\tn{H}_{m,\delta} |^2
   + \frac{\kappa(\vt)}{\vt} |\Grad \vt |^2 \right),
\]
which is still an {\it equality}\/ at this level.

Integrating \eqref{c3}, we deduce that
\bFormula{c3a}
{\rm ess} \sup_{t \in (0,T)} \ \| \log(\vt) \|_{L^1(\Omega)} \leq c,\
\| \log(\vt) \|_{L^2(0,T; W^{1,2}(\Omega))}\leq c.
\eF
Finally, we need an estimate on the time derivative of $\vt$. To get this,
we go back to \eqref{p6} (or, equivalently, multiply \eqref{c3} by $\vt$)
which we rewrite in the form
\bFormula{c3b}
 \partial_t \Big( \vt + M_\delta (\vt)  G(\tn{Q}) \Big)
  + \vu \cdot \Grad \Big( \vt + M_\delta (\vt)  G(\tn{Q}) \Big)
   - \Div \left( \kappa(\vt) \Grad \vt \right)
\eF
\[
 = \big( M_\delta(\vt) - U_\delta'(\vt) \vt \big)
  \big (\partial_t G(\tn{Q}) + \vu \cdot \Grad  G(\tn{Q}) \big)
\]
\[
  +  \left( \frac{\mu(\vt)}{2} \left| \Grad \vu + \Grad^t \vu \right|^2
  + \delta |\Grad \vu |^r + [\Gamma(\vt)]_\ep |\tn{H}_{m,\delta} |^2 \right),
\]
where we have set
\[
  M_\delta(\vt):= \int_1^{\vt} U_\delta''(r) r \ \dr.
\]
Then, testing \eqref{c3b} by a test function $\phi\in H^3(\Omega)$
and using the previous estimates (\ref{p21}--\ref{p23b2}) and
\eqref{c1}, and the fact that, by (\ref{i3}--\ref{i3bis}),
$M_\delta(\vt) \sim \vt^{1/2}$ for large $\vt$,
it is easy to check that
\bFormula{c3c}
  \big\| \partial_t \big( \vt + M_\delta (\vt)  G(\tn{Q}) \big) \big\|_{L^1(0,T;H^{-3}(\Omega))} \le c,
\eF
uniformly with respect to all approximation parameters. Since
$ G$ and $M_\delta$ are smooth and nonnegative and $\partial_t\tn{Q}$
is already estimated in \eqref{p23d}, we finally obtain
\bFormula{c3z}
 \| \partial_t \vt \|_{L^1(0,T;H^{-3}(\Omega))} \le c.
\eF

%%%%%%%%%%%%%%%%%%%%%%%%%%%%%%%5%%%%%%%%%%%%%%%%%%%%%%%%%%%%%%%%%%%%%%%%%%%%%%%%%%%%%%%%%%%%%%%%%%5

\subsection{The limit $m \to \infty$}

It is convenient to start with the limit for $m \to \infty$, with the other parameters
$N$, $\ep$, and $\delta$ fixed. In such a way, the velocity field remains regular at
this stage, with facilitates the limit passage considerably. As a result of this step,
we obtain that $f(\tn{Q})$ lies in $L^\infty(0,T; L^1(\Omega))$ in the
limit. This yields, in particular, uniform boundedness of $\tn{Q}$ by
a constant independent of the parameters $N$, $\ep$, and $\delta$.

We start by recalling that, in accordance with (\ref{p8}),
\bFormula{c4}
  {\rm ess} \sup_{t \in (0,T)} \intO{ f_m (\tn{Q}_m) } \leq c,
\eF
whereas, by virtue of the standard identity (see for instance \cite{ROCK})
\[
  \mathcal{L} \left[ \frac{\partial f_m (\tn{Q}_m)} {\partial \tn{Q}} \right] : \tn{Q}_m = f_m
   (\tn{Q}_m) + f^*_m \left( \mathcal{L} \left[ \frac{\partial f_m (\tn{Q}_m)} {\partial
   \tn{Q}} \right] \right),
\]
combined with (\ref{p22}) and \eqref{p23b}, we infer that
\bFormula{c5}
  \| f_m (\tn{Q}_m) \|_{L^{5/3}((0,T) \times \Omega)} \leq c.
\eF
We are ready to perform the limit for $m \to \infty$. To this end, we denote
by $\{ \vu_m , \tn{Q}_m, \vt_m \}_{m > 0}$
the family of approximate solutions. First, by virtue of \eqref{c4} and \eqref{c5},
we have
\[
  f_m(\tn{Q}_m) \to \Ov{f(\tn{Q})} \ \mbox{weakly-(*) in}\ L^\infty(0,T; L^1(\Omega))
    \ \mbox{and}\ \mbox{weakly in}\ L^{5/3}((0,T) \times \Omega),
\]
at least for suitable subsequences.

Next, using (\ref{p10}), (\ref{p23}), \eqref{p23b}, \eqref{p23d}, and the Aubin-Lions lemma (cf. \cite{aubin, lions}), we get
\[
\tn{Q}_m (t, \cdot) \to \tn{Q}(t, \cdot) \ \mbox{in, say,} \ L^2(\Omega;\TS) \ \mbox{for a.e.} \ t \in [0,T],
\]
\[
\Grad \tn{Q}_m \to \Grad \tn{Q} \ \mbox{(strongly) in}\
L^q((0,T) \times \Omega; R^{27}) \ \mbox{for all}\ q \in [1,10/3),
\]
and
\[
\Delta \tn{Q}_m \to \Delta \tn{Q} \ \mbox{weakly in}\ L^2((0,T) \times \Omega; \TS).
\]

Our ultimate goal is to show that
\bFormula{c7}
\tn{Q}(t, x) \in \mathcal{D}[ f ] \ \mbox{for a.a.} \ (t,x) \in (0,T) \times \Omega,
\eF
yielding, in particular, the desired conclusion
\[
\Ov{ f(\tn{Q}) } = f(\tn{Q}).
\]

Using (\ref{c5}), the monotonicity of the sequence $f_m$, the pointwise convergence of $\{ \tn{Q}_m \}_{m > 0}$,
and Fatou's lemma, we obtain
\[
{\rm ess} \sup_{t \in (0,T)} \intO{ f_{m_0} (\tn{Q}) } \leq c \ \mbox{for any fixed} \ m_0 , \ \mbox{with}\
c \ \mbox{independent of}\ m_0;
\]
whence, by the Levi theorem applied for $m_0 \to \infty$, we conclude that
\[
{\rm ess} \sup_{t \in (0,T)} \intO{ f (\tn{Q}) } \leq c,
\]
yielding (\ref{c7}).

Let us notice that, as a consequence of \eqref{c7}, we have in particular
\bFormula{c7b}
  \| \tn{Q} \|_{L^\infty((0,T) \times \Omega; \TS)} \le c.
\eF
Recalling \eqref{p23} and applying the standard Gagliardo-Nirenberg inequality (cf.~\cite[p.~125]{nier}),
we then also have
\bFormula{p23bis}
  \| \Grad\tn{Q} \|_{L^4((0,T) \times \Omega; R^{27})} \leq c,
\eF
uniformly in $N$, $\delta$ and $\ep$. As a consequence, we can also improve
(\ref{p23c}--\ref{p23d}) as follows:
\bFormula{p23c2}
  \big\| \vu \cdot \Grad \tn{Q} - \SSS(\Grad \vu, \tn{Q}) \big\|_{L^{20/11}((0,T)\times \Omega; \TS)} \leq c,
\eF
\bFormula{p23d2}
  \| \partial_t \tn{Q} \|_{L^{20/11}((0,T)\times \Omega; \TS)} \leq c.
\eF
%

%%%%%%%%%%%%%%%%%%%%%%%%%%%%%%%%%%%%%%%%%%%%%%%%%%%%%%%%%%%%%%%%%%%%%%%%%%%%%%%%%%%%%%%%%%%%%%%%%%%%

\subsubsection{Strong convergence of the temperature and the limit system }

Using (\ref{c1}--\ref{c2}), \eqref{c3z}, and the Aubin-Lions lemma (cf. \cite{aubin, lions}),
we deduce that
\[
  \vt_m \to \vt \ \mbox{in}\ L^q((0,T) \times \Omega) \ \mbox{for a certain}\ q > 1.
\]
Then, the uniform estimates derived above permit to take the limit $m\to \infty$ in
the momentum equation to obtain
\bFormula{c9}
\frac{{\rm d}}{{\rm d}t} \intO{ \vu \cdot \vc{v} }
\eF
\[
= \intO{ [\vu]_\delta \otimes \vu : \Grad \vc{v} } - \delta \intO{
|\Grad \vu |^{r-2} \Grad \vu : \Grad \vc{v} } - \intO{ \mu(\vt) (
\Grad \vu + \Grad^t \vu) : \Grad \vc{v} }
\]
\[
+ \intO{ \Big( \Grad \tn{Q} \odot \Grad \tn{Q} \Big) : \Grad \vc{v}
}
\]
\[
- \intO{ \left\{  2 \xi \left( \tn{H}_\delta : \tn{Q} \right) \left(
\tn{Q} + \frac{1}{3} \tn{I} \right) - \xi \left[ \tn{H}_\delta \left(
\tn{Q}  + \frac{1}{3} \tn{I} \right) + \left( \tn{Q}  + \frac{1}{3}
\tn{I} \right) \tn{H}_\delta \right] + \left( \tn{Q} \tn{H}_\delta - \tn{H}_\delta
\tn{Q} \right) \right\} : \Grad \vc{v} },
\]
for any $\vc{v} \in X_N$, where
\bFormula{pHdelta}
 \tn{H}_{\delta} = \Delta \tn{Q} - \mathcal{L} \left [
  \frac{\partial f (\tn{Q})} {\partial \tn{Q}} \right] + U_\delta(\vt) \mathcal{L} \left [
   \frac{\partial G (\tn{Q})} {\partial \tn{Q}} \right].
\eF

Analogously, we can take the limit $m\to\infty$ in the director
equation, obtaining
\bFormula{c10}
 \partial_t \tn{Q} + (\vu \cdot \Grad) \tn{Q} - \tn{S} (\Grad \vu, \tn{Q} ) =
  [\Gamma (\vt )]_\ep \tn{H}_\delta.
\eF
Actually, we remark that $\vu$ is still a smooth vector field in
the limit. On the other hand, we cannot prove strong convergence for $\tn{H}$ and
$\Grad\vt$. For this reason, the heat equation {\it cannot pass}\/ to the limit $m\to\infty$
in the form \eqref{p6}. For this reason, it is convenient to replace it
by a {\it partial form}\/ of the energy balance at this stage.

To do this, for fixed $m$ we multiply \eqref{p4} by $-\tn{H}_{m,\delta}$ and
sum the result to \eqref{p6}. Note that we {\it do not sum}\/ the energy contribution
coming from the momentum equations at this level. A number of integrations by parts
similar to those performed in Section~\ref{sub:bounds} then permit to deduce
\bFormula{c11m}
  \partial_t \left( \frac{1}{2} | \Grad \tn{Q}_m |^2 + f(\tn{Q}_m) - \Big( U_\delta (\vt_m)
   - \vt U'_\delta(\vt_m) \Big) G(\tn{Q}_m) + \vt_m \right)
\eF
\[
 + \Div \bigg( \vu_m \Big( \frac{1}{2} | \Grad \tn{Q}_m |^2 + f(\tn{Q}_m)
  - \Big( U_\delta (\vt_m) - \vt U'_\delta(\vt_m) \Big) G(\tn{Q}_m) + \vt_m \Big) \bigg)
\]
\[
  - \Div \Big( \kappa(\vt_m) \Grad \vt_m \Big) - \Div \Big( \Grad \tn{Q}_m : \tn{S}(\Grad \vu_m, \tn{Q}_m) \Big)
   - \Div  \Big( [\Gamma(\vt_m)]_{\ep} \Grad \tn{Q}_m : \tn{H}_{m,\delta} \Big)
\]
\[
  + \tn{S} (\Grad \vu_m, \tn{Q}_m ) : \tn{H}_{m,\delta}
   + \Big( \Grad \tn{Q}_m \odot \Grad \tn{Q}_m \Big) : \Grad \vu_m
  = \frac{\mu(\vt_m)}{2} \left| \Grad \vu_m + \Grad^t \vu_m \right|^2
   + \delta |\Grad \vu_m |^r.
\]
Now it is possible to take the limit $m\to\infty$ in the above relation.
Actually  the worst terms are the quadratic ones on the \rhs. However,
they do pass to the limit since the velocity still takes values in the finite
dimensional space $X_N$ at this level (and, hence, it is a smooth
function). We then get
\bFormula{c11}
  \partial_t \left( \frac{1}{2} | \Grad \tn{Q} |^2 + f(\tn{Q}) - \Big( U_\delta (\vt)
   - \vt U'_\delta(\vt) \Big) G(\tn{Q}) + \vt \right)
\eF
\[
 + \Div\bigg( \vu \Big( \frac{1}{2} | \Grad \tn{Q} |^2 + f(\tn{Q})
  - \Big( U_\delta (\vt) - \vt U'_\delta(\vt) \Big) G(\tn{Q}) + \vt \Big) \bigg)
\]
\[
  - \Div \Big( \kappa(\vt) \Grad \vt \Big) - \Div \Big( \Grad \tn{Q} : \tn{S}(\Grad \vu, \tn{Q}) \Big) - \Div
   \Big( [\Gamma(\vt)]_{\ep} \Grad \tn{Q} : \tn{H}_\delta \Big)
\]
\[
  + \tn{S} (\Grad \vu, \tn{Q} ) : \tn{H} + \Big( \Grad \tn{Q} \odot \Grad \tn{Q} \Big) : \Grad \vu
   = \frac{\mu(\vt)}{2} \left| \Grad \vu + \Grad^t \vu \right|^2 + \delta |\Grad \vu |^r.
\]
Due to the lack of strong convergence for $\tn{H}$ and $\Grad \vt$,
also the $m$-limit of the entropy equation \eqref{c3} has
to be written in the form of an inequality. Indeed, a
standard semicontinuity argument yields
\bFormula{c12}
  \partial_t \Big( \log(\vt) + U'_\delta (\vt) G(\tn{Q}) \Big)
   + \vu \cdot \Grad \Big( \log(\vt) + U'_\delta (\vt) G(\tn{Q}) \Big) -
    \Div \left( \frac{\kappa(\vt)}{\vt} \Grad \vt \right)
\eF
\[
  \ge \frac{1}{\vt} \left( \frac{\mu(\vt)}{2} \left| \Grad \vu + \Grad^t \vu \right|^2
   + \delta |\Grad \vu |^r + [\Gamma(\vt)]_\ep |\tn{H} |^2
   + \frac{\kappa(\vt)}{\vt} |\Grad \vt |^2 \right).
\]

%%%%%%%%%%%%%%%%%%%%%%%%%%%%%%%%%%%%%%%%%%%%%%%%%%%%%%%%%%%%%%%%%%%%%%%%%%%%%%%%%%%%%%%%%%%%%%5%%%%%%%%%%%%%%%

\subsection{The limit $N \to \infty$}
\label{sec:N}

Our next goal is to let the Galerkin parameter $N \to \infty$.
In this step, the most difficult points consist
in taking the limit of (\ref{c9}), in order to let it
converge to the Navier-Stokes system,
and of the internal energy equation (\ref{c11}).
To achieve these limits, we need to prove a strong
convergence for $\Grad \vu_N$. With this purpose,
we denote $\{ \vu_N, \tn{Q}_N, \vt_N \}_{N > 0}$
the family of approximate solutions and recall that
\bFormula{c13}
 \| \tn{Q}_N \|_{L^\infty((0,T) \times \Omega; \TS)} \leq c,
\eF
with $c$ independent of $N$, $\ep$, and $\delta$.
Then, we take $\vc{v}=\vu_N$ in \eqref{c9} and
correspondingly test \eqref{c10} by
$\tn{H}_N$. Repeating the usual cancellations, we then arrive at
\bFormula{semi-1}
  \partial_t \intO{ \left( \frac{1}{2} | \vu_N |^2
   + \frac12 | \Grad \tn{Q}_N |^2 + f( \tn{Q}_N) \right) }
   - \intO{ U_\delta(\vt_N) \mathcal{L}
      \left[\frac{\partial G (\tn{Q}_N)} {\partial \tn{Q}_N} \right] :
       \big( \partial_t \tn{Q}_N + \vu_N \cdot \Grad \tn{Q}_N \big) }
\eF
\[
  + \intO{ \left(\delta | \Grad \vu_N |^r
   + \frac{\mu(\vt_N) }{2} \left|\Grad \vu_N + \Grad^t \vu_N \right|^2
   + [\Gamma(\vt_N)]_\ep | \tn{H}_N |^2 \right) }
   = 0,
\]
where we have used the identities \eqref{mat}
%
%\bFormula{mat}
% - \tn{H} : \tn{S}(\Grad \vu_N, \tn{Q}_N)
%\eF
%\[
% = \left( \tn{Q}_N \tn{H} - \tn{H} \tn{Q}_N \right) : \Grad \vu _N
% + 2 \xi \left( \tn{H}: \tn{Q}_N \right) \left( \tn{Q}_N : \Grad \vu_N \right)
%- \xi \left[ \tn{H} \left( \tn{Q}_N  + \frac{1}{3} \tn{I} \right) +
%\left( \tn{Q}_N  + \frac{1}{3} \tn{I} \right) \tn{H} \right]: \Grad \vu_N,
%\]
%
%that holds for any symmetric matrix $\tn{H}$,
%
and
\bFormula{mat1}
  \intO{ ( \vu_N \cdot \Grad \QQ_N ) : \Delta \QQ_N
   + ( \Grad \QQ_N \odot \Grad \QQ_N ) : \Grad \vu_N }
   = 0.
\eF
Next, we integrate \eqref{semi-1} in time and compute its supremum
limit as $N\to \infty$. We obtain
\bFormula{semi-2}
  \limsup_{N\to\infty} \int_0^T \intO{ \left(\delta | \Grad \vu_N |^r
   + \frac{\mu(\vt_N) }{2} \left|\Grad \vu_N + \Grad^t \vu_N \right|^2
   + [\Gamma(\vt_N)]_\ep | \tn{H}_N |^2 \right) } \ \dt
\eF
\[
   = \lim_{N\to \infty} \int_0^T
   \intO{ U_\delta(\vt_N) \mathcal{L}
      \left[\frac{\partial  G (\tn{Q}_N)} {\partial \tn{Q}_N} \right] :
       \big( \partial_t \tn{Q}_N + \vu_N \cdot \Grad \tn{Q}_N \big) } \ \dt
\]
\[
  - \liminf_{N\to\infty} \intO{ \left( \frac{1}{2} | \vu_N(T) |^2
   + \frac12 | \Grad \tn{Q}_N(T) |^2 + f( \tn{Q}_N ) \right) }
  + \intO{ \left( \frac{1}{2} | [\vu_0]_\delta |^2
   + \frac12 | \Grad [\tn{Q}_0]_\delta |^2 + f( \tn{Q}_0 ) \right) }
\]
\[
   \le \intO{ U_\delta(\vt) \mathcal{L}
      \left[\frac{\partial  G (\tn{Q})} {\partial \tn{Q}} \right] :
       \big( \partial_t \tn{Q} + \vu \cdot \Grad \tn{Q} \big) } \ \dt
\]
\[
  - \intO{ \left( \frac{1}{2} | \vu(T) |^2
   + \frac12 | \Grad \tn{Q}(T) |^2 + f( \tn{Q} ) \right) }
  + \intO{ \left( \frac{1}{2} | [\vu_0]_\delta |^2
   + \frac12 | \Grad [\tn{Q}_0]_\delta |^2 + f( \tn{Q}_0 ) \right) }.
\]
Actually, the existence of the limit in the second line above is ensured by
the pointwise convergence of $\vt_N$, by the boundedness of $U_\delta$
and $G$ together with their first derivatives, and by
estimates \eqref{c1} and (\ref{p23c2}--\ref{p23d2}).
Moreover, the liminf on the
third line above can be estimated by
using semicontinuity of norms with respect to weak convergence.

Now, we observe that it is possible to take the limit $N\to \infty$ in
\eqref{c10}. Analogously, we can write the limit of \eqref{c9};
however, the second integral on the right
hand side has to be temporarily written as
\bFormula{semi-2b}
  - \delta \intO{ \Ov {|\Grad \vu |^{r-2} \Grad \vu} : \Grad \vc{v} }
\eF
in the limit, since strong convergence of $\Grad \vu$
is not achieved yet.

Then, we take $\vc{v} = \vu$ in the $N$-limit of \eqref{c9}
and test the $N$-limit of \eqref{c10} by $\tn{H}$.
Note that this procedure is rigorous
(and gives rise to an equality) thanks to the fact that
we have that $\Grad \vu$ is bounded in $L^r$, where $r>3$,
uniformly in $N$ and, consequently, the coupling term
\[
  \vu \cdot \Grad \tn{Q} - \SSS(\Grad \vu, \tn{Q})
\]
in \eqref{c10} lies in $L^2$ even in the limit.
Hence, equation \eqref{c10} can still
be read as a relation in $L^2$ and use of the $L^2$-test
function $\tn{H}$ is consequently permitted (note that this will
no longer be true in the limit $\delta\to 0$).
This is exactly the reason why the $r$-Laplacean regularization
has been added in the momentum equation.

This procedure permits to achieve, in the
limit $N\to \infty$, the analogue of \eqref{semi-1}.
Comparing with \eqref{semi-2} we then obtain
\bFormula{semi-3}
  \limsup_{N\to\infty} \int_0^T \intO{ \left(\delta | \Grad \vu_N |^r
   + \frac{\mu(\vt_N) }{2} \left|\Grad \vu_N + \Grad^t \vu_N \right|^2
   + [\Gamma(\vt_N)]_\ep | \tn{H}_N |^2 \right) } \ \dt
\eF
\[
  \le \int_0^T \intO{ \left(\delta \Ov{|\Grad \vu |^r \Grad \vu} \cdot \Grad \vu
   + \frac{\mu(\vt) }{2} \left|\Grad \vu + \Grad^t \vu \right|^2
   + [\Gamma(\vt)]_\ep | \tn{H} |^2 \right) } \ \dt,
\]
whence a simple monotonicity argument permits to
deduce the strong convergence
\bFormula{semi-3b}
  \Grad \vu_N \to \Grad\vu \ \text{in}\ L^r((0,T)\times \Omega; R^9).
\eF
Using this relation together with \eqref{p23bis}, \eqref{p23d2},
and the Aubin-Lions lemma (cf. \cite{aubin, lions}), we can also take the limit
$N\to \infty$ of \eqref{c11} (and in particular of the
quadratic terms on the \rhs).

Moreover, \eqref{semi-3b} permits to identify the limit of
${|\Grad \vu |^{r-2} \Grad \vu}$ in the momentum equation
(cf.~\eqref{semi-2b}). In particular, the Navier-Stokes
system \eqref{c9} passes to the desired limit as $N\to \infty$.

Next, we notice that the (infimum) limit $N\to \infty$
can be taken in the entropy inequality
\eqref{c12} as in the previous section. Thus,
to complete the passage to the limit w.r.t.~$N\to \infty$
it is sufficient to recover the total energy balance
(cf.~\eqref{i24}). With this aim, we can notice that,
in the limit $N\to \infty$, \eqref{c9} is no longer a relation
in the finite-dimensional space $X_N$, but can rather be interpreted
as a true PDE, (at least) in the distributional sense.
In other words, we can rewrite it in the form
\bFormula{i22x}
 \partial_t \vu + \Div ([\vu]_\delta \otimes \vu) =
   \Div \sigma,
\eF
where
\bFormula{i25x}
  \sigma = \mu(\vt) \left( \Grad \vu + \Grad^t \vu \right)
   + \delta | \Grad \vu |^{r-2} \Grad \vu
  - p \tn{I}
\eF
\[
  + 2 \xi \left( \tn{H}: \tn{Q} \right) \left( \tn{Q} + \frac{1}{3} \tn{I} \right)
  - \xi \left[ \tn{H} \left( \tn{Q}  + \frac{1}{3} \tn{I} \right)
  +  \left( \tn{Q}  + \frac{1}{3} \tn{I} \right) \tn{H} \right]
  + \left( \tn{Q} \tn{H} - \tn{H} \tn{Q} \right) - \Grad \tn{Q} \odot \Grad \tn{Q}.
\]
Then, multiplying \eqref{i22x} by $\vu$ and adding the result to
\eqref{c11}, standard integrations by parts permit to get
the {\it total energy balance}\/ \eqref{i24} (where of course,
at this level, the stress $\sigma$ still contains the
regularizing contribution $\delta | \Grad \vu |^{r-2} \Grad \vu$).
\bRemark{onreg}
 In order to achieve the total energy balance \eqref{i27}
 it is crucial to remove first the Galerkin approximation (otherwise
 the kinetic energy contribution is projected on the finite-dimensional
 space $X_N$). For this reason we need that, after taking the
 limit $N\to \infty$, we still have sufficient regularity to
 use $\vu$ as a test function in the momentum equation
 in order to get a kinetic energy {\it equality}.
 This regularity is properly provided by the additional term
 $\delta | \Grad \vu |^{r-2} \Grad \vu$.
\eR

%%%%%%%%%%%%%%%%%%%%%%%%%%%%%%%%%%%%%%%%%%%%%%%%%%%%%%%%%%%%%%%%%%%%%%%%%%%%%%%%%%%%%%%%%%%%%%5%%%%%%%%%%%%%%%

\subsection{The limits $\varepsilon\to 0$ and $\delta\to 0$}

After taking the limit $N\to \infty$, we still have to remove the
regularizations coming from the parameters $\ep$ (appearing in the
coefficient $[\Gamma(\vt)]_\ep$ in \eqref{p4}), and $\delta$
(appearing in the mollification of the initial data,
in the function $U_\delta$ in the approximation of $\HH$,
and in the regularizing terms $[\vu]_\delta \otimes \vu$ and
$\delta | \Grad \vu |^{r-2} \Grad \vu$ in the momentum
equation).

Since the total energy balance \eqref{i27} is already achieved
and no quadratic terms are present in its \rhs, letting $\ep\to 0$
and subsequently $\delta\to 0$ does not give rise to any additional
difficulty. Actually, most of the argument can be carried
out just by adapting the procedure, based on the uniform a-priori
bounds of Section~\ref{a}, used before to let
$m,N\to \infty$. The only point which requires some
additional care is letting
$\delta\to 0$ in the momentum equation. However, this procedure
is absolutely analogous to the argument outlined in
\cite[Sec.~5.2]{FFRS} to which we refer the reader for
details. We just note that we need here the restriction $r<10/3$
(cf.~Section~\ref{pp}) on the exponent of the additional viscosity term
(up to this point we only used that $r>3$).

\bibliographystyle{alpha}

\end{document}